\documentclass[11pt]{amsart}
\usepackage{geometry}                
\usepackage[parfill]{parskip}    
\usepackage{graphicx}
\usepackage{amssymb}
\usepackage{epstopdf}
\usepackage[all, 2cell]{xy}
\usepackage{enumitem}
\usepackage{eufrak}
\usepackage{ebproof}
\usepackage[linktocpage]{hyperref}
\usepackage{amsmath}

\hypersetup{
   colorlinks=true, 
   linktoc=all,     
   linkcolor=blue,  
}
\DeclareRobustCommand{\subtitle}[1]{\\#1}

\title{Kolmogorov's Calculus of Problems and Its Legacy}

\author{Andrei Rodin}

\date{\today}

\begin{document}
\maketitle
\tableofcontents

\renewcommand{\subtitle}[1]{}

\begin{abstract}
Kolmogorov's Calculus of Problems is an interpretation of Heyting's intuitionistic propositional calculus published by A.N. Kolmogorov in 1932. Unlike Heyting's intended interpretation of this calculus, Kolmogorov's interpretation does not comply with the philosophical principles of Mathematical Intuitionism. This philosophical difference between Kolmogorov and Heyting implies different treatments of problems and propositions: while in Heyting's view the difference between problems and propositions is merely linguistic, Kolmogorov keeps the two concepts apart and does not apply his calculus to propositions.
I stress differences between Kolmogorov's and Heyting's interpretations and show how the two interpretations diverged during their development. In this context I reconstruct Kolmogorov's philosophical views on mathematics and analyse his original take on the Hilbert-Brouwer controversy. Finally, I overview some later works motivated by Kolmogorov's Calculus of Problems and propose a justification of Kolmogorov's distinction between problems and propositions in terms of Univalent Mathematics.      
\end{abstract}

\newpage
\section{Introduction}
 The \emph{Calculus of Problems} was proposed by Andrei N. Kolmogorov in his paper titled \emph{On the Interpretation of Intuitionistic Logic} written and published originally in German in 1932 \cite{Kolmogoroff:1932}, English translation \cite{Kolmogorov:1998}. The paper contains no technical results but develops an argument according to which the symbolic calculus  published by Arent Heyting earlier in 1930 under the name of \emph{Intuitionistic Logic} should be interpreted 
 \footnote{The term ``interpretation'' in such contexts should not be understood in the model-theoretic sense when one interprets a formal theory in terms of another theory; the term ``explanation'' would be here more appropriate \cite{Sundholm:1983}, \cite{vanAtten:2022}. Having this in mind I nevertheless stick to the more common usage in this paper.}
 in terms of problems rather than sentences (propositions) and called accordingly
  \footnote{In 1925 Kolmogorov independently published his own formalisation of intuitionistic logic \cite{Kolmogorov:1925}, English translation \cite{Kolmogorov:1991a}, which in today's terms can be described as the \emph{minimal} fragment of the predicate Intuitionistic logic. In his 1932 paper Kolmogorov does not refer to his 1925 paper but use the more recent Heyting's formalisation \cite{Heyting:1930}, English translation \cite{Heyting:1998}. Kolmogorov's 1925 paper was published in Russian and was unknown to Heyting in 1930. In an undated letter to Heyting sent after 1933 \cite[p.15-16]{Troelstra:1990}, English translation \cite[p.91]{Kolmogorov:1988}, Kolmogorov mentions his 1925 paper but it is unlikely that Heyting later could read it. For detailed overviews of both Kolmogorov's logical papers see \cite{Troelstra:1990}, \cite{Uspenskii&Plisko:1991}, \cite{Coquand:2007}. }.
   
In his comprehensive monograph \cite{Heyting:1934} Heyting acknowledges the difference between his and Kolmogorov's interpretation of the intuitionistic logic but in later writings he describes his and Kolmogorov's interpretations of intuitionistic logic as essentially the same \cite{Heyting:1958}. Following Heyting Troelstra later coined the name and the concept of the  BHK-interpretation of the intuitionistic logic (so called after the names of Brouwer, Heyting and Kolmogorov), with the idea that Heyting and Kolmogorov  ``have an equal share''  in this unified interpretation  \cite[p.7]{Troelstra:1990} 
 \footnote{The acronim ``BHK'' first appears in press in 1977 but at this point letter ``K'' stands for Georg Kreisel but not for Kolmogorov. Here I refer to the standard version of the BHK-interpretation as it is presented in \cite{Troelstra:1990} and many other contemporary and later publications. See \cite{vanAtten:2022} for further historical details.}. 
 
It is not immediately clear whether or not Kolmogorov, in his turn, was agree to consider his Calculus of Problems (CP) as essentially the same as Heyting's Intuitionistics Logic in its mature form, say, as presented (and explained) in \cite{Heyting:1934}. In what follows I analyse the available evidence, and reconstruct Kolmogorov's position. I show that unlike Heyting's Intuitionistic Logic (IL) Kolmogorov's CP was not intended to belong to the \emph{intuitionistic mathematics}, i.e., in mathematics developed on philosophical principles of Mathematical Intuitionism. Kolmogorov's idea, instead, was to extend the existing logical analysis of mathematical reasoning, which was focused on mathematical propositions or assertions, with a new logical theory dealing with mathematical \emph{problems}. I show, in what follows, that this Kolmogorov's project was consistent with his philosophical views on mathematics, and that his view on IL remained stable throughout his long career. In the nutshell the difference between Kolmogorov's and Heyting's interpretations of IL is this. The difference is both philosophical and technical. Philosophically, Kolmogorov rejects Mathematical Intuitionism as a viable philosophy of mathematics, while Heyting follows in Brouwer's steps. Technically, Kolmogorov insists on the distinction between problems and propositions while Heyting (in \cite{Heyting:1934} and later on) treats the two notions as essentially the same. 

I conclude that Heyting's 1934 view on Kolmogorov's interpretation of IL is correct but Heyting's 1958 judgement on the same matter is not accurate. The issue is important because the vast majority of later historians of intuitionsitic logic and mathematics followed Heyting's 1958 judgement, which caused a misrepresentation of Kolmogorov's contribution. Having said that, I don't  challenge the notion of BHK-interpretation, which is perfectly coherent theoretically whatever is its name, but try to fix a historical misconception to which this popular notion easily leads.

One may wonder whether Kolmogorov's unorthodox interpretation of IL has a theoretical value. In Section 4 I argue that it does. For this end I point to two recent developments, namely to combined calculus of problems and propositions QHC due to Sergei Melikhov \cite{Melikhov:2022a}, \cite{Melikhov:2022b} and to Homotopy type theory (HoTT) \cite{UF:2013}, where the distinction between problems and propositions suggested by Kolmogorov shows up. The two developments are very different in their character. QHC is an explicit attempt do realise Kolmogorov's project suggested in his \cite{Kolmogoroff:1932} and in the related sources. By contrast, the influence of Kolmogorov's CP on HoTT is only indirect, namely, via the Inuitionisitc Type theory due to Per Martin-L\"of \cite{Martin-Lof:1984} (MLTT). As we shall see in Section 4 the differences between Kolmogorov's and Heyting's interpretations of IL, which I'm trying to elucidate here, could hardly play a role in the development of MLTT; following Heyting Martin-L\"of  treats here propositions and problems  as essentially the same notion. Most certainly, Kolmogorov's CP didn't play any specific role in the discovery of homotopical interpretation of MLTT by Awodey and Voevodsky, which gave rise to HoTT.  It is moreover striking to observe that the distinction between problems and propositions suggested by Kolmogorov is supported by HoTT for a purely mathematical reason.


\section{Kolmogorov's Calculus of Problems and Its Philosophical Background}
The presentation of Kolmogorov's CP given in this Section is intended to be self-contained but it highlights first of all epistemological aspects of Kolmogorov's work and thus doesn't aim at completeness. In particular, I'm trying to stress the aspects of Kolmogorov's work, which make his approach divergent from Heyting's. For other presentations and analyses of CP I refer the interested reader to   \cite{Troelstra:1990}, \cite{Uspenskii&Plisko:1991}, \cite{Coquand:2007}, and \cite{Melikhov:2017}, which provide different perspectives on the same subject. In what follows I refer to these works repeatedly discussing some specific points.     

\subsection{Overview}
Kolmogorov's paper \cite{Kolmogoroff:1932} is divided into two parts. In the first part Kolmogorov introduces CP and stresses that it is not constrained by the principles of Mathematical Intuitionism; in the second part of the paper some intuitionistic concepts including the intuitionistic notion of negation proposed by Brouwer are critically discussed.

The purpose of CP is explained in the first part as follows:
 
\begin{quote}
Along with the development of theoretical logic, which systematises the schemes of proofs of theoretical truths, it is also possible to systematise the schemes of solutions of problems.   \cite[p.328]{Kolmogorov:1998}, translation corrected
\end{quote}  

More than a half century later during the preparation of his Collected Works Kolmogorov made the following comment on his \cite{Kolmogoroff:1932} where he reiterated the same thesis:
 
\begin{quote}
Paper ``On the interpretation of intuitionistic logic'' was written with the hope that the logic of solutions of problems would later become a regular part of courses on logic. It was intended to construct a unified logical apparatus dealing with objects of two types \textemdash\ propositions and problems. \cite[p.452]{Kolmogorov:1991} \footnote{The comment first appears in press in Russian in 1985 as a part of Russian edition of Kolmogorov's Collected Works.}.
\end{quote}

The \emph{reducibility} of a problem to another problem in CP is analogous to implication in propositional logic:
\begin{quote}
If we can reduce the solution of problem $b$ to the solution of problem $a$, and the solution of problem $c$ to the solution of problem $b$, then the solution of $c$ can also be reduced to the solution of $a$. (ibid)
\end{quote} 

The concepts of problem and solution of a problem Kolmogorov does not define but illustrates with examples. Here are two examples: 
\begin{quote}
\begin{enumerate}
\item To find four whole numbers $x, y, z, n$ for which the relations
$$x^{n}+y^{n}=z^{n},\quad  n>2$$
hold.
\item To prove the falsity of  Fermat's [Last] Theorem\footnote{Fermat Last Theorem says that there exist no quadruple of whole numbers satisfying the conditions of Problem 1. This theorem was conjectured by Pierre Fermat in 1637, remained open at the time of writing of Kolmogorov's paper and was proved in 1995 by Andrew Wiles.}. \cite[p.329]{Kolmogorov:1998}
\end{enumerate}
\end{quote}

Kolmogorov remarks that (2) reduces to (1), i.e., that a solution of (1) solves (2), but not the other way round since (2) can be proved by deriving a contradiction from the statement of the theorem without finding the quadruple of numbers satisfying (1). He stresses once again that pointing to this difference between problems (1) and (2) ``does not yet constitute a particular intuitionistic claim'' \cite[p.329]{Kolmogorov:1998}, i.e. does not depend on epistemological principles of Mathematical Intuitionism. At the same time, the corresponding \emph{propositions}

 \begin{enumerate}[label=(\roman*)]
\item There \underline{exist} four whole numbers $x, y, z, n$ for which the relations
$$x^{n}+y^{n}=z^{n},\quad  n>2$$
hold.
\item Fermat's [Last] Theorem is false.
\end{enumerate}
are obviously equivalent, see \cite[p.333, note 5]{Kolmogorov:1998}. 

\subsection{Negated Problems and Unsolvable Problems} 

By negation of problem $a$, in symbols $\neg a$, Kolmogorov understands ``the problem 'to obtain a contradiction provided that the solution of $a$ is given'' \cite[p.329]{Kolmogorov:1998}, which corresponds to the standard intuitionistic negation used by Heyting. This explanation of negation is followed by an interesting footnote:
\begin{quote}
Let us observe that $\neg a$ should not be understood as the problem ``to prove the unsolvability of [problem] $a$''. If one considers in general ``the unsolvability of $a$'' as a well-defined concept, than one only obtains the theorem that from $\neg a$ follows the unsolvability of $a$, but not the converse. For example, if it were proven that the well-ordering of the continuum surpasses our abilities, one could still not claim that a contradiction follows from the existence of such a well-ordering \cite[p.333, note 5]{Kolmogorov:1998}.  
\end{quote}  

The above remark leaves a room for different interpretations. Coquand interprets it in the sense that the existence of well-ordering of the continuum (or the existence of some other mathematical construction) can be not provable in some theory (say, ZF) but still be consistent with this theory (which is indeed the case of ZF in view of the independence of the Axiom of Choice relatively to this theory) \cite[p.31]{Coquand:2007}. According to Melikhov, in this remark Kolmogorov ``conflates implications [i.e, reductions] between problems with implications between propositions'', which makes the remark ambiguous \cite[p.24]{Melikhov:2017}. My analysis by and large agrees with Melikhov's but I'll try to present it here in simpler terms without using Melikhov's logical framework and references to independence results. The main problematic point of this remark, in my view, concerns the notion of \emph{solution} (of a given problem), which Kolmogorov leaves undefined. For a simpler example of problem let us consider the trisection of a given angle into equal parts with ruler and compass and call it problem $T$. Now by \emph{solution} of $T$ one may understand either
\begin{enumerate}[label=(\alph*)]
\item the two wanted lines trisecting the given angle or 
\item a general method of constructing such lines with ruler and compass for any given angle.
\end{enumerate}

Clearly only (b) but not (a) corresponds to the usual notion of what is an expected solution of $T$. But observe that (a) is on par with the intended solution of problem (1) above, where no specific constructive means are specified (as the wanted quadruple of numbers).  As we know it today via the Galois theory, $T$ is \emph{unsolvable} in the sense that, provably, no solution of $T$ of type (b) exists. In other words, (a hypothesis of) the existence of (b) leads to contradiction. But the existence of solution in the sense (a) is intuitively obvious and does not lead to contradiction. Now, if we understand ``$T$ is unsolvable'' as usual (i.e. as above) but define $\neg T$ as ``the existence of solution (a) leads to contradiction'' we get the situation described in the above Kolmogorov's remark: $\neg T$ does \underline{not} hold but $T$ is nevertheless unsolvable. This example strongly suggest to redefine $\neg T$ as ``the existence of solution (b) leads to contradiction''. In this case $\neg T$ and ``$T$ is unsolvable'' become equivalent. 

But there is a further and deeper reason Kolmogorov might see this latter solution unsatisfactory. It has to do with general problems of intuitionistic (constructive) negation thoroughly discussed in the second part of Kolmogorov's paper  \cite[S. 64-65]{Kolmogoroff:1932}. A systematic analysis of this problematic issue is out of the scope of the present paper, and I illustrate here the difficulty only with the chosen example. A basic idea of the intuitionsitic mathematics is to avoid reasoning about mathematical objects without referring to well-defined procedures that allow one to build (construct) such objects as mental entities. Even if Kolmogorov didn't want to commit himself to the epistemological principles of intuitionistic mathematics, he wanted to apply the same approach in the limited domain of mathematical problems \textemdash\ without giving up the classical reasoning with mathematical propositions. Yet, in order to define the notion of negated problem he relies on the idea of \emph{existence} of its solution, which is not supported by any concrete construction. The talk of existence of solutions of type (b) rather than of type (a) apparently doesn't alleviate the difficulty but only opens an infinite regress: instead of a non-constructive talk about lines or numbers one gets here a talk about possible constructive procedures (say, about possible constructions with ruler and compass), which is equally non-constructive. By 1932 David Hilbert was pursuing the idea of using a similar move for reducing various non-constructive mathematical theories involving abstract non-intuitive ``ideal'' objects to an elementary theory (under the names of Proof theory or Meta-mathematics) dealing only with intuitive finite ``real'' syntactic constructions, which encode the ``ideal'' theories \cite{Hilbert&Bernays:1934-1939}. In Hilbert's view the reduction of ``ideal'' mathematical objects and theories to their ``real'' (finitary) counterparts was a clear epistemic gain. But I can see no sign that Kolmogorov in \cite{Kolmogoroff:1932} followed a similar agenda. His footnote in question is rather an evidence to the contrary. This is why I have some doubts that using a modern metamathematical machinery is appropriate for interpreting this Kolmogorov's remark. This gives theoretically interesting outcomes but it may be simply reading too much into it. In the context of Kolmogorov's elementary examples of problems, a non-constructive talk about possible constructive procedures has no clear epistemic advantage with respect to non-constructive talks about mathematical objects themselves.

\subsection{Postulates and Axioms according to Kolmogorov}        

For saving space I will not present here Kolmogorov's interpretation of IL in terms of problem-solving systematically. It is as anyone familiar with the standard informal BHK-semantics would expect. I shall comment here only on some interesting stylistic aspects of Kolmogorov's presentation. His list of the axioms (tautologies) of IL is preceded by the following commentary:

 \begin{quote}
We \emph{postulate} that we have already solved the following two groups of problems $\dots$. The further presentation addresses only a reader who has already solved all these problems \cite[p.330]{Kolmogorov:1998}.
 \end{quote} 

\emph{First}, it is remarkable that Kolmogorov uses here the term ``postulate'' (``Wir nehmen als Postulaten an'' in the original German \cite[S.61]{Kolmogoroff:1932}) in the the same exact sense as Euclid in his \emph{Elements}: recall that at least Euclid's basic Postulates 1-3, which define the rules of constructions ``with ruler and compass'' do not have a propositional form but are elementary problems or operations (such as ``to construct a straight line between two given non-identical points'') that help to solve further problems and prove further theorems (with various ``auxiliary'' geometrical constructions) \cite[ch. 2]{Rodin:2014}. 

\emph{Second}, the solutions of these elementary problems are not simply taken for granted as solvable but claimed to be effectively solved by the author. The same is required from the reader (more on this point will be said shortly). In a footnote Kolmogorov compares CP with a propositional calculus in this respect as follows: 

 \begin{quote}
In the case of propositional calculus, one must first convince oneself of the correctness of of the axioms, if one want to determine the correctness of the consequences \cite[p.334, note 9]{Kolmogorov:1998}.
 \end{quote}

``Solving a postulate'' or ``proving an axiom'' (or ``convincing oneself'' that a given axiom is ``correct'', if one prefers) is, obviously, a different kind of epistemic task than solving a regular problem or proving a regular theorem. The former unlike the latter does not involve using known problems and theorems for solving a new problem or proving a new theorem. But however these elementary epistemic acts are further explained, the above quote makes it clear that Kolmogorov sticks to a traditional Frege-style notion of axiom (and an Euclid-style notion of postulate), and doesn't follow Hilbert's new ideas about the axiomatic method and axiomatic theories, of which by 1932 Kolmogorov is very well aware
\footnote{
Kolmogorov provides a critical discussion about Hilbert's axiomatic method in his philosophical paper \emph{Contemporary Debates on the Nature of Mathematics} first published in Russian in 1929 \cite{Kolmogorov:1929}, see English translation \cite[p. 380, 383, 385]{Kolmogorov:2006}
}
. 

 \emph{Third}, Kolmogorov's appeal to his own and the reader's individual mental experience is quite remarkable. Plausibly, here Kolmogorov is influenced by his reading of Brouwer. Earlier in the same paper Kolmogorov explains his view on the individual mental experience more precisely:

\begin{quote}
The fact that \emph{I} solved a problem is a purely subjective fact that in itself has as yet no general interest. However, the logical and mathematical problems possess the special property \emph{of the general validity} [der Allgemeing\'ultigkeit] \emph{of their solutions}: If I have solved a logical or mathematical problem, then I can present this solution in a way that is intelligible to all and it is \emph{necessary} that it be recognised as correct solution although this necessity has to a certain extent an ideal character, for it presupposes a sufficient intelligence on the part of the listener. \cite[p.330]{Kolmogorov:1998}
\end{quote} 

The above clarification is at odds with Brouwer's view on mathematics as an individual mental activity.  
\footnote{Cf. Brouwer's thesis: ``Intuitionism $\dots$ highlights the existence of pure mathematics independent of language'' \cite[p. 50]{Brouwer:1998}
}
Kolmogorov takes seriously not only the subjective mental experience of an individual mathematician but also the communicational aspect of mathematical problem-solving that allows a mathematician to share their relevant experiences with others and even, idealistically, with ``all''. Thus asking the reader to ``solve'' all the ``postulates'' of CL (i.e., the axioms of IL) the Kolmogorov wants to make sure that his individual mental experience is shared with the reader properly, and generates an appropriate individual experience on the side of the reader. It is not immediately clear what the author can do in order to make sure that his individual experience is shared properly. More detailed verbal translations and explanations of each Kolmogorov's ``postulate'' (which are missing in the paper but can be made either by a commentator or by the reader themselves) can be certainly helpful. Notice also, that the author has no means for checking the results of his efforts unless he gets a feedback from his readers. But however this mechanism works, it is, after all, an empirical sociological fact that mathematical and logical understanding \emph{is} sharable to a certain extent, and the case in point is not an exception.   

In a following footnote Kolmogorov extends  the same argument, \emph{mutatis mutandis}, to `` the proof of theoretical propositions'' but remarks that while every proven proposition is conventionally called ``correct'' [richtig]
\footnote{
According to today's standard, the talk of ``correct'' proposition in the sense suggested by Kolmogorov is very informal.  Pointing to a mathematical proposition that has a generally accepted proof, today's working mathematician would rather call it ``true'' while a logician would use words ``provable'' or ``proved'' relativising these properties to a relevant foundational framework. At the same time the talk of ``correct proof'' can make part of a professional exchange both in the general mathematical community and in the logical community.     
}
, this latter term has no counterpart in the case of a problem and its solution \cite[p.334, note 8]{Kolmogorov:1998}. Indeed, the English adjective ``correct'' can be applied to solution of a given problem but not to the problem itself (unless it is used in the sense of ``well-posed'' problem but this is not what Kolmogorov means here). The same remarks apply to German adjective ``richtig''. So the missing term for problems should be synonymous to ``solved'' or ``solvable'' (the latter term to be understood in the sense of an easily realised possibility). My guess is that Kolmogorov's motivation behind this terminological remark is his decision to give to the notions of problem and proposition equal rights in logic and mathematics, and  by all means avoid their conflation.

\subsection{Law of Excluded Middle}
Let us now see how Kolmogorov's explains why the Law of Excluded Middle (LEM, in symbols $\vdash a \vee \neg a$) does not universally hold for problems. Under Kolmogorov's interpretation LEM is understood as the problem 

``to solve a given problem $a$ or to prove that $a$ has no solution (without premises)''. 

In the special case when $a$ has a form ``to prove proposition $p$'', LEM is interpreted as problem

 ``to prove $p$ or lead $p$ to contradiction''.
 
Since $a$ is an \emph{arbitrary} problem we don't have here a room for ambiguity between the notion of solution as a method and that as a wanted mathematical object. The only candidate for a solution of LEM is a general \emph{method} $M$ allowing one to solve (in the positive sense of obtaining a solution or in the negative sense of proving that no solution exists) \emph{any problem whatsoever}. In the propositional case the solution amounts to a general method that allows one either to prove or disprove (by reducing to absurdity) any mathematical proposition. Now, according to Kolmogorov's epistemic standard explained above, in order to claim LEM (under the problem-solving interpretation) to be solved, the author should possess  $M$ in some form and be able to communicate it to others. On that point Kolmogorov remarks that 

 \begin{quote}
If our reader does not consider himself to be omniscient, he will probably determine that [LEM] cannot be found on the list of problems solved by him. \cite[p.332]{Kolmogorov:1998}
 \end{quote}     
 
The above argument against LEM, in Kolmogorov's view, does not apply to \emph{propositions}.

Kolmogorov's strategy to deal with Brouwer's general critique of LEM is made clear in the following remark made in Kolmogorov's philosophical paper of 1929 
 \cite{Kolmogorov:1929}, English translation \cite{Kolmogorov:2006}  

\begin{quote}
The law of excluded middle according to Brouwer could not be applied only to a certain kind of judgements, in which a theoretical statement is closely connected with construction of the object of the statement. Therefore, we may assume that Brouwer's ideas do not contradict the traditional logic, which has never before dealt with such judgements. \cite[p 385]{Kolmogorov:2006}, translation corrected
\end{quote}  
 
Most plausibly by ``certain kind of judgements'' Kolmogorov  means here judgements involving Brouwer's so-called ``weak counter-examples''  to the classical reasoning. The weak counter-examples are mathematical  constructions that involve some open problems (say, the Goldbach Conjecture)  \cite{vanAtten:2020} like ``real numbers that do not have a decimal expansion''  \cite[33-34]{Brouwer:1998a}. Unlike Brouwer and Heyting Kolmogorov does not think that these examples demonstrate the need to revise all of the ``ordinary'' mathematics and disqualify many instances of classical mathematical reasoning using LEM.  Kolmogorov does not treat Brouwer's examples explicitly but a strategy of doing this suggested by Kolmogorov's CP is to distinguish between \emph{problem} $P$ and \emph{proposition} ``problem $P$ has a solution'', and apply LEM universally in case of propositions but not in case of problems \cite{Melikhov:2022a}. Then one is in a position to argue, contra Brouwer, that every real number, including those numbers which according to Brouwer at the present stage of mathematical knowledge do not have decimal expansions (since some problems involved in the construction of extension remain open), in fact, \emph{do} have such  extensions \textemdash\ quite independently of one's knowledge of these extensions and of one's capacity to effectively expand the given number indefinitely. It goes without saying that this obvious remark does not disqualify Brouwer's neat analysis of decimal extensions of different kinds of real numbers provided in his {Brouwer:1998a}. Brouwer's talk of ``decimal numbers without decimal extension'' is very appealing and provocative but it is not absolutely necessary for appreciating his ingenious mathematical results.  

In order to develop properly Kolmogorov's line of defence of classical mathematical reasoning with LEM against Brouwer's objections, one needs a theory supporting a rigorous mathematical and logical distinction between problems and propositions rather than simply relying on conventional meaning of these words borrowed from the mathematical practice. In Section 4 we'll overview a recent work by Sergei Melikhov who provides a foundation of such a theory. In Section 3 we'll see how Kolmogorov and Heyting treated and discussed  this important distinction before Heyting wholly abandoned it in his \cite{Heyting:1934}.       

\subsection{Kolmogorov and Brouwer-Hilbert Controversy}
I conclude this current Section with some remarks on the philosophical background behind Kolmogorov's CP.  Kolmogorov's repeated insistence on the fact that CP does not depend on philosophical principles of Mathematical Intuitionism can be easily misread as a claim of philosophical neutrality. One can get an impression that Kolmogorov unlike Brouwer and Heyting did not pursue any specific philosophical agenda but simply tried to interpret IL in common mathematical terms of his time where the notions of problem and its solution obviously belonged
\footnote{
Compare this interesting comparison of Kolmogorov's 1932 paper with his earlier 1925 Russian paper where Kolmogorov proposes a formalisation of intuitionistic logic and constructs a version of double-negation translation of classical formulas into intuitionistic formulas (see \cite{Coquand:2007} for an overview):   

\begin{quote}
In \cite{Kolmogorov:1925} an embedding operation is constructed which makes it possible to give an intuitionistic interpretation to the major part of classical mathematics, while, in a sense, the paper \cite{Kolmogoroff:1932} is devoted to the solution of the inverse problem of interpreting intuitionistic logic within the framework of ordinary mathematical notions, irrespective of the philosophical and methodological principles of intuitionism. \cite[p.459]{Uspenskii&Plisko:1991}, the references to Kolmogorov's papers are adjusted to the bibliography of the present paper  
\end{quote} 

}. 
But the thesis of Kolmogorov's alleged philosophical neutrality is straightforwardly refuted by the available evidence. In his philosophical paper \emph{Contemporary Debates on the Nature of Mathematics} first published in Russian in 1929  \cite{Kolmogorov:1929} (English translation  \cite{Kolmogorov:2006}) Kolmogorov provides a  critical review of Brouwer-Hilbert controversy and the related foundational debates. He compares the contemporary situation in the foundations of mathematics with that of the mid 19th century, and judges that while in the 19th century foundations of Differential and Integral Calculus could be fixed by purely mathematical and logical means, the ``solutions of the present controversies should be found beyond the scope of mathematics'' \cite[p. 380]{Kolmogorov:2006}. 
Discussing the problem of non-constructive existence in mathematics Kolmogorov criticises a popular philosophically neutral approach to this problem:

\begin{quote}
The most common issue from this difficulty among mathematicians that avoid philosophy is limitation of the domain of ``existence''. [$\dots$]. This position \textemdash\ though the most placid one \textemdash\ suffers from unprincipledness, which is expressed most evidently in the fact that bounds of what each mathematician is ``ready to admit'' depend on his personal interests. \cite[p. 383]{Kolmogorov:2006}   
\end{quote}

Thus ``avoiding philosophy'' could hardly be Kolmogorov's motivation behind his distancing from the Mathematical Intuitionism in his 1932 paper  \cite{Kolmogoroff:1932}. What was then Kolmogorov's philosophical position in foundations of mathematics? Analysing Hilbert's and Brouwer's arguments in \cite{Kolmogorov:1929}, Kolmogorov identifies a common problem in these two lines of thought as follows:

\begin{quote}
The emergence of these extreme points of view [i.e., of Hilbert's Formalism and Brouwer's Intuitionism] could be explained by the fact that the combination of the two [namely, formal and constructive] aspects  of the set-theoretic mathematics caused great difficulties and even contradictions. The common source of these difficulties is the following. Mathematicians have been used to treating numbers, functions, and sets as if they were things of the real world, similar to material things. The very preference of the word ``thing'' (Ding) to the word ``object'' (Gegenstand) [as in Hilbert's \emph{Grundlagen der Geometrie} of 1899 \cite{Hilbert:1899}] is rather typical in this sense.  \cite[p. 382]{Kolmogorov:2006}, translation corrected
\footnote{I corrected the translation by replacing ``objects'' to ``things'' (Russian ``veshch'' in the third sentence of the quote. This is important in view of Kolmogorov's point made in the last sentence of the same quote. Kolmogorov's distinction between things and objects is discussed below in the main text.} 
\end{quote}

This Kolmogorov's argument goes along with a remark made by Kolmogorov much later, in 1975, during his public lecture in Moscow State University, which is reported by V.A. Uspenskii.  A listener asked Kolmogorov what the foundations of mathematics were about. Kolmogorov (according to Uspensii):
\begin{quote}
Mathematics studies objects that do not quite really exist. The foundations of mathematics take care of the transition from the experience to these abstractions, in order to not get confused. \cite[p. 301]{Uspenskii:2006}
\end{quote} 

The above argument is original and to the best of my knowledge has no evident analogue in the vast secondary literature on the Brouwer-Hilbert controversy. Let me try to explain and expend it. In his \emph{Foundations of Geometry} \cite{Hilbert:1899} Hilbert famously begins with suggesting the reader to think of ``systems of things'' [Systeme von Dingen] of unspecified nature, which satisfy certain conditions that Hilbert calls ``axioms''. Under the intended interpretation the ``things''  become usual Euclidean points, straight lines and plane. But they can be also represented (modelled) by some different mathematical objects (in particular, by numbers and arithmetical constructions) or avoid any specific representation and remain abstract ``thought-things'' [Gedankendinge], which ``exist'' merely in virtue of the fact that one can think of them consistently \cite[p.43]{Rodin:2014}. This is the core of Hilbert's axiomatic method, which became very influential in the 20th century mathematics, logic, and philosophy \textemdash\ but which, as we have seen above, Kolmogorov does not accept. The possibility to ``make up'' mathematical objects as explained above, impose on them some consistent sets of relational properties (axioms), and prove theorems on this basis is the first target of Kolmogorov's critique. In his view, which can be described as ``traditional'' with respect to Hilbert's non-traditional approach, axioms are justified with properties of corresponding mathematical objects but not the other way round. Kolmogorov:

\begin{quote} 
In order to give sense to an abstract [Hilbert-style axiomatic] theory, the existence of at least one system of objects and relations satisfying the proposed axioms is required. When systems with a finite number of elements are considered, the question is easily solved, because such a system could be materially presented. \cite[p. 381]{Kolmogorov:2006} 
\end{quote}

So in Kolmogorov's view, the existence of a model of a given axiomatic theory  not only implies a relative consistency property for this theory but also and foremost allows one to ``make sense'' of this theory and, in particular (in view of his remarks on mathematical axioms in \cite{Kolmogoroff:1932} quoted above), allows one to establish that the axioms of the given theory are \emph{true} (with respect to the given model). Kolmogorov well understands Hilbert's notions of non-interpreted axiom and non-interpreted theory but he doesn't grant to them the same epistemic significance as Hilbert. Notice also Kolmogorov's reference to the ``material'' representation in the above quote. As we shall shortly see, it is important for his wanted conception of mathematical object. I would like also to stress that Kolmogorov's conservative position with respect to Hilbert's axiomatic method is combined with his full awareness about all major developments in his contemporary ``abstract'' mathematics including developments in Set theory in which he contributed himself. 
\footnote{
Kolmogorov's rejection of Hilbert's Formalism is shared by his friend and colleague of the time Alexandr Khintichin who in 1926 published a philosophical paper \emph{Ideas of intuitionism and the struggle for a subject matter in contemporary mathematics} \cite{Khinchin:1926}, English translation in \cite{Verburgt&Hoppe-Kondrikova:2006}. The ``subject matter'' in the title is English translation of Russian ``predmet'', which can be also translated into English as ``object'' and into German as ``Gegenstand'' \textemdash\ notice the same German word used by Kolmogorov in his critique of Hilbert in the above quote. Khinchin's critique of Hilbert's Formalism is similar to Kolmogorov's but unlike Kolmogorov Khinchin does not provide in this paper also critical arguments against Brouwer's Intuitionism.   
}.  

Hilbert's roundabout way to treat mathematical contents involving  infinities via finitary syntactic constructions, which they represent, Kolmogorov describes as a ``brilliant art'' (or ``utmost skill'' as translated in  \cite[p. 386]{Kolmogorov:2006}) which, however, ``gives no explanation on how mathematics has existed up to the moment and how mathematicians could understand each other while uttering senseless statements on the infinite'' (ib.). Kolmogorov is more sympathetic to Brouwer's project, which, in his view, does ``not fear facing the problem and promise to reveal the nature of the infinite'' (ib.). Yet, Kolmogorov believes that Brouwer's proposed solution is erroneous:
\begin{quote} 
[O]ne may doubt that the intuition and the construction of new forms from positive integers would prove to be reliable in this case. In particular, Brouwer explores the continuum in the form of infinite sequences of positive integers $\dots$. However, historically the idea of the continuum has been created by idealization of a really observable continuous environment. Now it is hard to imagine how we could find in this a basis for the development of mathematical theory. However, only this could be a straight way to understand the nature of mathematical continuum.
\end{quote}  

Even the above quote concerns only Brouwer's theory of continuum, the quote makes it clear that Kolmogorov doesn't share basic tenets of Brouwer's Intuitionism and  thinks about the nature of mathematical objects in a different way. Kolmogorov states this explicitly in his Preface to the Russian translation of \cite{Heyting:1934}, which appeared in 1936  \cite{Kolmogorov:1936}, see Appendix for the full English translation: 

 \begin{quote}
We cannot agree with the intuitionists when they claim that mathematical objects are products of the constructive activity of our spirit. For us, mathematical objects are abstractions from existing forms of reality, which is independent from our spirit.  (Appendix)
\end{quote}

Let me now turn to the core of Kolmogorov's argument according to which a \emph{common} mistake of Hilbert and Brouwer is that both these mathematicians uncritically follow the usual professional mathematical parlance and think about mathematical objects  ``as if they were things of the real world, similar to material things''  \cite[p. 382]{Kolmogorov:2006}. In a similar (but more specific) context Kolmogorov talks about ``an overrealistic attitude towards “things” with which the mathematical theory deals''  \cite[p. 383]{Kolmogorov:2006}

The charge of naive realism (Platonism) about mathematical objects made against Hilbert and against Brouwer may appear very surprising. Indeed, Hilbert's notion of \emph{thought-thing} related to his version of the axiomatic method doesn't suggest a realistic and moreover naturalistic interpretation. On the contrary, Hilbert's notion of axiomatic theory helped many mathematicians in the 20th century to think  of mathematical theories as being wholly independent of theories of physics or any other natural science. 
\footnote
{To give just one example, in 1961 Marshall Stone explained the virtues of modern axiomatic method with the following strong statement:
\begin{quote}  
While several important changes have taken place since 1900 in our conception of mathematics or in our points of view concerning it, the one which truly involves a revolution in ideas is the discovery that mathematics is entirely independent of the physical world.'' \cite[p.716]{Stone:1961}. 
\end{quote}    
}. 
Brouwer's notion of mathematical object as mental construction does not suggest any direct analogy with material objects either. So what Kolmogorov could have in his mind when he claimed that both Hilbert and Brouwer were mislead by the false analogy between mathematical objects and material objects? A possible charitable reading of Kolmogorov's argument is suggested by his distinction between \emph{thing} (German \emph{Ding}) and \emph{objects} (German \emph{Gegenstand}) made in this context
\footnote{The German equivalents of words ``thing'' and ``object'' are provided in the Russian original text of \cite{Kolmogorov:1929} by Kolmogorov.}.

As Kolmogorov makes it explicit in his paper, the ``things'' refer to Hilbert's \emph{Foundations of Geometry} \cite{Hilbert:1899} as explained above. But what does he mean by ``objects'', German ``Gegenst\"ande'', in this context? Even if one cannot rule out the possibility that this Kolmogorov's remark is merely linguistic, it strongly suggests that Kolmogorov is knowledgable of some contemporary German philosophical and psychological literature that treats the notion of  Gegenstand as a key epistemological concept. 

A likely (albeit not only possible) Kolmogorov's source is the \emph{Theory of Objects} ( Gegenstandstheorie) developed by Alexius Meinong  during the first decade of the 20th century, see \cite{Meinong:1904} (English translation \cite{Meinong:1960}), \cite{Meinong:1907}. In the recent literature Meinong is often misrepresented as a champion of ontological inflation populating ontology with all sorts of possible and impossible objects. Hence the popular pejorative use of the expression``Meinong's jungles'' first used by W. Kneale \cite{Jacquette:1996}. This is a sheer misreading of Meinong, however, since his main epistemological idea was in a sense the opposite: in order to be an object of study this object needs not to exist or even be logically possible \cite{Jacquette:1996}, \cite{Routley&Routley:1973}. Mathematical objects served to Meinong as important examples of such non-existing theoretical objects (and in fact, he even considered the pure mathematics as a proper part of his Gegenstandstheorie \cite{Meinong:1907}). A Meinongian reading of Kolmogorov's objection to Hilbert is as follows. Whatever ontological status Hilbert gives to his \emph{thought-things} and other mathematical objects, he gives them \emph{some} ontological status, i.e. assumes that they in some sense \emph{exist}. This is where an analogy with material objects plays a role: Hilbert applies to mathematical objects the same logical principles, which in different contexts apply to reasoning about material objects or about any other kind of \emph{entities}. But in Kolmogorov's (and in Meinong's) view mathematical objects are not entities at all, and thus require a wholly different logical treatment. What is at stakes here is the basic logical notion of \emph{individual} used by Hilbert in foundations of mathematics but not this or that view on the ontological status of mathematical objects
\footnote{Cf. Routley\&Routley:
The real perniciousness of Platonism is that, by his existence assumptions, the Platonist is enable to transfer to non-entities a logical structure suitable only to entities and evolved in that case. \cite[p.246]{Routley&Routley:1973}.}.
In Kolmogorov's  view, the false analogy between mathematical \emph{objects} and material ``things'' aka entities is responsible for Russell paradox  and difficulties about the Axiom of Choice \cite[p. 383]{Kolmogorov:2006}. His strategy, as I understand it, is to develop an appropriate theory of mathematical objects and only on that basis develop foundations of mathematics including its logical foundation. Whether or not Meinong's works indeed motivated Kolmogorov's stress on the concept of Gegenstand, Kolmogorov's remarks on foundations of mathematics including the aforementioned Uspenskii's evidence of 1975 are perfectly consistent with Meinong's theory. 

The fact that Kolmogorov makes the same charge (of relying on the false analogy between mathematical objects and material entities) against Brouwer's Intuitionism may appear even more surprising. By 1929 (the date of publication of Kolmogorov's philosophical paper  \cite{Kolmogorov:1929}) Brouwer already developed a detailed philosophical account of mathematics (as a part of his general philosophical view on life) where natural numbers and other admissible mathematical objects were explained away in terms of ``primordial intuition'' of time, human freedom and a number of other philosophical and moral concepts, some of which were first introduced by Brouwer himself \cite{vanStigt:1990}. This philosophical view on mathematics is highly controversial and can be objected on various grounds but the charge of naive realism  about mathematical objects made against this view doesn't make any sense. 

When Kolmogorov was preparing his philosophical paper \cite{Kolmogorov:1929} for publication he could hardly be aware about the full philosophical background of Brouwer's revisionary program in the foundations of mathematics because the first Brouwer's publication in German containing such details appears the same year (publication 1929A in the list of Brouwer's publication found in  \cite[ch.1] {vanStigt:1990}
\footnote{The same reference system is used in the Collected Works by Brouwer edited by A. Heting \cite{Brouwer:1975}. The same reference system is used in the present paper also for other Brouwer's works. For English translation of Brouwer's 1929A see \cite{Brouwer:1998}. 
}.
Earlier Brouwer's publications of the same purely philosophical character (1905 \emph{Life, Art and Mysticism} and 1919B \emph{Mathematics, Truth, Reality}) were all in Dutch, and could be hardly accessible to Kolmogorov, who was proficient only in German, French and English in addition to his native Russian. Like the rest of the international mathematical community of the time, Kolmogorov learned of Brouwer's program in foundations of mathematics via his papers such 1921A  \emph{Does Every Real Number Have a Decimal Expansion?} which combined philosophical insights with a mathematical ingenuity but didn't contain a systematic philosophical argumentation. Should we on this basis simply discharge the aforementioned Kolmogorov's objection to Brouwer as irrelevant and based on a sheer misunderstanding of Brouwer's position?       

I don't think so. In fact, in  \cite{Kolmogorov:1929} Kolmogorov makes it very clear what he has here in his mind. He objects here against the notion of mathematical object as  ``construction $\dots$ based on positive integers or on some other resource of elementary objects'' \cite[p. 382]{Kolmogorov:2006}. The objection is that \emph{this} constructivist notion is based on a false analogy between mathematical objects and material entities. While in Hilbert's case the relevant material entities are any physical objects, in Brouwer's case the relevant kind of material object is a technical artefact. True, this objection is not specific to Brouwer's Mathematical Intuitionism and applies more straightforwardly to Kronecker's dictum ``God created the natural numbers; all the rest is the work of man'' and to any other variety of Mathematical Constructivism. Brouwer's notion of mathematical construction is, of course, far from being philosophically naive and it is certainly not naturalistic. Yet, when it comes to a mathematical work, the Mathematical Intuitionism comes down to the same idea of using natural numbers (and perhaps some other kinds of distinguished elementary objects) and some well-specified procedures applied to the elementary objects as the only legitimate way of introducing further mathematical objects. Kolmogorov criticises this core idea of Mathematical Constructivism quite independently of its specific philosophical underpinning that one can find in Brouwer's writings. Even if Kolmogorov's argument concerning the false analogy with material constructions would hardly make sense to Brouwer, it could be appealing to many contemporary mathematicians interested in foundations of their discipline who didn't share Brouwer's philosophical ideas and convictions.

Kolmogorov's research in foundations of mathematics and in its philosophy like Hilbert's research in these fields is led by the idea to justify what Kolmogorov perceived as the best mathematical practices of the time. This included works in Set theory and related abstract areas of modern mathematics, that is, all of the ``Cantor Paradise''. This is why along with Hilbert Kolmogorov rejects Brouwer's revisionary program that imposes new sever restrictions on admissible ways of mathematical reasoning. But Hilbert's formalist strategy  of grounding mathematics on finitary syntactic structures, as I have already explained,  is quite unsatisfactory in Kolmogorov's eyes either. Along with Brouwer Kolmogorov is looking for an appropriate contentful notion of mathematical object. He does not accept Brouwer's strategy of solving this problem but he thinks that Brouwer's and Heyting's works shed a new light on an area of mathematical practice, which in the contemporary discussion so far remained outside of logical analysis, namely, the problem-solving. Hence Kolmogorov's idea to develop on the basis of  Heyting's work a new Calculus of Problems as a complement to the existing Calculus of Propositions. In 1936 Kolmogorov describes this double approach as follows:

 \begin{quote}
[C]onstructive solutions of problems are as much important in mathematics as the pure proofs of theoretical sentences. This constructive aspect of mathematics does not conceal for us its other and more fundamental aspect, namely, its epistemic aspect. But the laws of mathematical construction discovered by Brouwer and systematised by Heyting under the appearance of new intuitionistic logic, so understood, preserve for us their fundamental significance. see Appendix
\end{quote} 

It should be stressed that Kolmogorov's writings available to the date do not contain anything like a systematic theory of mathematical objects
\footnote{Kolmogorov's archive is presently a private property, and at the time of writing it has not been fully studied. This leaves a hope that some relevant documents can become available in the future.}. 
Yet, the textual evidences quoted above allows one to understand, by and large, Kolmogorov's position in the debates on foundations of mathematics during the first half of the 20th century. I summarise the above analysis in form of a list of desiderata for the hypothetical theory of mathematical objects (theory K), which Kolmogorov, in my understanding, aimed at:   

\begin{enumerate} 
\item according to K, mathematical objects are \emph{fictions} in the strong sense that a mathematical work with these objects involves strictly \emph{no} ontological commitment;  
\item yet, according to K, mathematical objects provide the mathematical thought with a genuine semantic \emph{content} aka a subject-matter;  
\item K supports a notion of \emph{mathematical truth} (rather than dispenses with truth in mathematics like the current versions of mathematical fictionalism \cite{Balaguer:2018};
\item K supports a notion of \emph{mathematical objectivity} and accounts for the social aspects of mathematics;  
\item K supports a theory of \emph{mathematical abstraction} and idealisation which is compatible with Scientific Realism and strongly links mathematics with the world of Physics and other natural sciences, including Cognitive Science;
 \item K explains away the conventional talk of mathematical existence;
\item K accounts both for the constructive (in particular, computational) aspect of mathematics (including problem-solving) and its non-constructive propositional aspect, and explains how the two are related.  
\item K supports the established mathematical practices, including the current research practices in modern abstract mathematics, and does not imply a need of their deep normative revision (as in the case of Mathematical Intuitionism).   
\end{enumerate} 

No theory satisfying all the above desiderata is known by the date. Noticeably, (1)-(5) are supported by Meinong's Theory of Objects albeit this theory lacks any sense of mathematical rigour.  To develop an account of mathematical objects in line of the above desiderata is, in my view, a challenging and interesting philosophical project, particularly, in view of new developments in Cognitive Science and in Computer Science.

\section{Heyting and Kolmogorov on Problems, Proofs and Propositions}
\subsection{Kolmogorov and Intuitionistic Mathematics}
Heyting's paper \cite{Heyting:1930} (English translation  \cite{Heyting:1998b}) used by Kolmogorov for developing CP in his \cite{Kolmogoroff:1932}, issues from an earlier unpublished essay, presently unavailable, written in response to a prize question  proposed in 1928 by G. Mannoury via the Dutch Mathematical Society. The question was  to codify with a ``formal mathematical system $\dots$ regularities $\dots$ which Brouwer uses to give expression to his mathematical intuition''  \cite[p.4]{Troelstra:1990}. Accordingly, Heyting begins his \cite{Heyting:1930}  with and explanation of the idea of ``intuitionistic mathematics'', which implies an important reservation with respect to his following attempt to codify the intuitionistic reasoning with a logical calculus:

\begin{quote}
Intuitionistic mathematics [Intuitionistische Mathematik] is a mental activity [Denkt\"atigkieit], and for it every language, including the formalistic one, is only a tool for communication. It is in principle impossible to set up a system of formulas that would be equivalent to intuitionistic mathematics \cite[p. 311]{Heyting:1998b}. 
\end{quote}

Instead of exploring further the notion of \emph{intuitionistic mathematics}, which Heyting develops in Brouwer's steps, let me stress here that the very idea of developing a special kind of mathematics called ``intuitionistic'' or otherwise is wholly alien to Kolmogorov. Kolmogorov refers repeatedly to the ``intuitionstic point of view'' in mathematics \cite{Kolmogoroff:1937}, to the ``intuitionistic epistemological  assumptions'' \cite{Kolmogoroff:1932}, to Brouwer's ``intuitionism'' (as an epistemological position) \cite{Kolmogorov:1929} and to the ``intuitionists'' (the people adhering to this position like Heyting) \cite{Kolmogorov:1929}, \cite{Kolmogorov:1936} but to the best of my knowledge he never talks about the ``intuitionstic mathematics'' as a special \emph{kind} of mathematics. Kolmogorov certainly acknowledges the existence of different philosophical views on and different foundational programs in mathematics but unlike Heyting he doesn't accept the idea that a particular philosophical view or a particular foundational program can or should produce anything like an autonomous school of mathematical thought. 

The independence of Kolmogorov's CP from the epistomological tenets of mathematical intuitionism has been not only officially stated in Kolmogorov's \cite{Kolmogoroff:1932} but also fully acknowledged by Heyting in 1934 in his \cite{Heyting:1934}. Heyting points here to the independence of CP of Mathematical Intuitionism as the major difference between Kolmogorov's interpretation of IL and his own interpretation. After presenting his intended interpretation of IL Heyting remarks:

 \begin{quote}
Kolmogorov developed an akin [verwandten] idea which, however, goes beyond the former idea [i.e., Heyting's intended interpretation] since it provides Heyting's calculus with a meaning that does not depend on the intuitionistic assumptions [intuitionistischen Voraussetzungen] \cite[14]{Heyting:1934}, my translation from German.  
\end{quote} 

In 1958, however, Heyting described Kolmgorov's interpretation of 1932 and his own contemporary interpretation as essentially the same:
 \begin{quote}
The older interpretations by Kolmogoroff (as a calculus of problems) and Heyting (as a calculus of intended construction) were substantially equivalent. \cite[p. 107]{Heyting:1958}
\end{quote} 

Why Heyting changed his view on Kolmogorov's CP over the years? I cannot point to any specific reason but a general reason seems to be clear. During the time span between 1934 and 1958 Heyting's contributed a very significant effort into the project of developing the \emph{intuitionisitic mathematics} as a special school of mathematical thought. Main results of this work are summarised in Heyting's 1956 book \cite{Heyting:1956}, which includes chapters on (the intuitionsitc versions of) Number Theory, Real Analysis, Algebra, Geometry, Measure Theory and Integration Theory, elements of Functional Analysis (Hilbert Spaces), and, finally, Logic
\footnote{The fact that the chapter on Logic is placed in the end of this presentation of the intuitionistic mathematics reflects the intuitionistic view on the role of logic in mathematics.}. This shows that for Heyting the talk of \emph{intuitionstic mathematics} was not just a convenient linguistic expression but a genuine research project, which he pursued during all his professional career.   

In his paper of 1958 addressed to a philosophical audience Heyting, in addition to explaining the philosophical background of Mathematical Intuitionism, develops a narrative about the genesis of the intuitionistic school, which covers Brouwer's foundational works and major later developments. Clearly, in this context  Heyting mentions Kolmogorov only as a contributor to the development of the intuitionistic mathematics \textemdash\ disregarding the fact that Kolmogorov himself hardly ever intended to contribute to this project. Now, as far as Kolmogorov's contribution to the intuitionistic mathematics is concerned, it is indeed the case that Kolmogorov's ideas behind CP, which we presented and analysed in the previous Section, simply don't matter: an intuitionist is in a position to disregard the fact that for Kolmogorov the problem-solving is just one aspect of mathematical thinking and reasoning, which should not be confused with theorem-proving (while for an intuitionist it amounts to the same, as we shall shortly see). In other words, an intuitionist is in a position to interpret Kolmogorov's CP in intuitionistic terms, and get back Heyting's IL with its original intuitionistic interpretation. This explains and justifies in a sense Heyting's 1958 remark about Kolmogorov quoted above. But this analysis also supports the claim that the above Heyting's 1934 comparison of his and Kolmogorov's interpretations of IL is  accurate while Heyting's 1958 remark on Kolmogorov is biased and inaccurate. Comparing Heyting's accounts of Kolmogorov's interpretation of IL of 1934 and of 1958 one cannot even say that Heyting changed his opinion about Kolmogorov's work over the passed years. In 1958 Heyting  boldly disregards rather than somehow re-evaluates the difference between the two interpretations of IL acknowledged by him in 1934 \textemdash\ apparently, simply because it no longer attracts his attention. It is a change of focus and of perspective rather than a change of opinion.  

Unfortunately, the same narrow view on Kolmogorov's CP is also characteristic for a significant number of later historians who equally aim at building the legacy of the intuitionistic school without making a sufficient effort (in Kolmogorov's case) to consider it the context of other contemporary developments. When van Dalen states in 1979 that         

\begin{quote}
[b]oth Heyting and Kolmogoroff's interpretation [Sic!] were fundamental in nature, i.e., they were intended as the ``true'' meaning of intuitionistic logic. Of the two, clearly Heyting’s interpretation is foundationally the more important one. Quoted after \cite[p.159]{Sundholm:1983} 
\end{quote}

he is quite right as far as the name of ``intuitionistic logic'' is understood at face value as a logic built upon the intuitionistic epistemological principles and codifying reasoning in the intuitionistic mathematics. But this remark wholly ignores the fact that in Kolmogorov's view the ``true'' meaning of IL lies outside of Mathematical Intuitionism and the intuitionistic mathematics.  

Sundholm sees the only difference between Kolmogorov's and Heyting's interpretations of IL in the fact the former interprets IL in terms of \emph{problems} while the latter interprets IL in terms of \emph{expectations} or \emph{intensions} \cite{Sundholm:1983}. Commenting on Heyting's 1958 account of Kolmogorov's interpretation Sundholm says that here ``Heyting himself came to see that his was not a different explanation from that of Kolmogoroff'' \cite[160]{Sundholm:1983}. I claim that this is a misreading. The difference between Kolmogorov's and Heyting's interpretations of IL, which has been rightfully stressed by Heyting in 1934 (the relevant quote is also found in Sundholm's paper), does not concern semantic differences between words ``problem'' and ``expectation'' or ``intention''. It concernes the fact that Kolmogorov's interpretation of IL does not comply with the intuitionistic viewpoint on mathematics. As for the word ``problem'', Heyting himself used it interchangeably with ``expectation'' as early as in 1930, see \cite[p. 234]{Heyting:1930a} and \cite[p. 307]{Heyting:1998} for English translation. In the absence of an evidence to the contrary, I assume that Heyting at this occasion used this ``problem interpretation of IL'' independently of Kolmogorov: 

\begin{quote}
A proposition $p$ like, for example, ``Euler constant [$C$] is rational'' expresses a \emph{problem} [(un probl\`eme)], or, better yet, a certain expectation [(une certaine attente)] (that of finding two integers $a$ and $b$ such that $C = \frac{a}{b}$), which can be [either] fulfilled (r\'ealis\'e) or disappointed (d\'e\c{c}ue). \cite[p. 307]{Heyting:1998}   
\end{quote} 

The term ``intention'', which has larger philosophical  connotations, first appears in the same context in a closely related Heyting's paper written in German of 1931 \cite{Heyting:1931}, English translation \cite{Heyting:1983}:

\begin{quote}
A mathematical proposition [(Aussage)] expresses a certain expectation [(Erwartung)].  For example, the proposition, ``Euler constant $C$ is rational'' expresses [(bedeutet)] the expectation that we could find two integers $a$ and $b$ such that $C = \frac{a}{b}$. Perhaps, the word ``intention''   [(Intention)], coined by the phenomenologists, expresses even better what is meant here. \cite[p. 58]{Heyting:1983}  
\end{quote} 

So it did \emph{not} take to Heyting nearly 30 years, as Sundholm suggests, to figure it out that the interpretations of IL in terms of problems and in terms of expectations or intentions amount to the same. This was quite clear to Heyting already in 1930, very likely before Heyting read a draft of Kolmogorov's 1932 paper \cite{Kolmogoroff:1932}
\footnote{In his Second letter to Heyting \cite[p. 90]{Kolmogorov:1988}  Kolmogorov thanks Heyting for having read an unpublished draft of his \cite{Kolmogoroff:1932} (see also \cite[p. 7]{Troelstra:1990}). In the same letter Kolmogorov refers to Heyting's paper \cite{Heyting:1930}, and reveals some details of his planned future talk and the International Congress of Mathematicians in Z\"urich scheduled for 1932. Since the letter is undated, its content does not allow one to conclude when Kolmogorov sent the unpublished manuscript to Heyting. Notice that the available evidence does not rule out the possibility that Kolmogorov could read \cite{Heyting:1930a} before publishing  \cite{Kolmogoroff:1932}, and thus borrow the idea of problem-based interpretation of IL from Heyting (rather then the other way round). However, the fact that in  \cite{Kolmogoroff:1932} Kolmogorov refers only to  \cite{Heyting:1930} but not to \cite{Heyting:1930a} makes this hypothesis rather implausible.  Most likely, in my view, the two mathematicians came to this idea independently.}

Van Atten in his comprehensive article on the history of Intuitionistic Logic \cite{vanAtten:2022} similarly leaves out of his consideration the fact (stressed both by Kolmogorov in \cite{Kolmogoroff:1932} and by Heyting in \cite{Heyting:1934}) that Kolmogorov's interpretation of IL, unlike Heyting's interpretation, is not supposed to comply with the tenets of Mathematical Intuitionism \cite{vanAtten:2022}. Van Atten acknowledges the fact that Kolmogorov in  \cite{Kolmogoroff:1932}, unlike Heyting in 1934 \cite{Heyting:1934}, insists on distinguishing between problems and propositions.  But without considering  Kolmogorov's non-constructive philosophical motivations behind this distinction it appears to be contingent and merely lingustic. Further, van Atten refers to the aforementioned Heyting's 1958 remark according to which Kolmogorov's and Heyting's interpretations of IL are essentially the same \cite[p. 107]{Heyting:1958}, and then suggests an argument according to which ``[b]y 1937, Kolmogorov seems to have come to believe the same''. To support his claim van Atten refers to a short (10 lines) Kolmogorov's review of a published exchange  between Freudenthal and Heyting where Kolmogorov explains the reader the intuitionistic treatment of propositions using  words ``intention'' and ``problem'' interchangeably: 
\begin{quote}
As is [well] known, from the intuitionist point of view, a mathematical proof always consists in carrying out a construction. If a sentence is uttered hypothetically without proof, then it is rather only an intention or a problem to find a corresponding construction. (\cite{Kolmogoroff:1937}, my translation from German). 
\end{quote}

If I understand the argument correctly, van Atten along with Sundholm in his \cite{Sundholm:1983} assumes that the difference between Kolmogorov's and Heyting's interpretations of IL reduces to semantic differences between term ``problem'' used by Kolmogorov in his interpretation of IL, and terms ``intention'' and ``expectation'' used by Heyting for the same purpose.  On the basis of this assumption \emph{and} the above clause where Kolmogorov uses ``problem'' and ``intention'' interchangeably, van Atten concludes that at the time of writing (1937) Kolmogorov did not any longer see any difference between his interpretation of IL and Heyting's. In my opinion this argument is erroneous because the above assumption is. As we have just seen, words ``problem'' and ``expectation'' have been used interchangeably in a similar context by Heyting back in 1930 but this didn't prevent him later from acknowledging the difference between his and Kolmogorov's interpretations  in 1934 \cite{Heyting:1934}. Without the above assumption van Atten's argument is not conclusive: the fact that Kolmogorov 1937 used terms ``problem'' and ``intention'' interchangeably does not imply that at this point he changed his view on IL as it was described by Heyting in 1934.    

Leaving the philosophical background aside, one can formulate the core difference between the two interpretations of IL as follows: according to Kolmogorov, IL applies only to problems but not to propositions (and the distinction between problems and propositions is essential) while, according to Heyting (at his mature stage of  \cite{Heyting:1934}), IL also applies to propositions because between propositions and problems there is no difference other than merely linguistic.  
In the above quote Kolmogorov makes it clear that he presents the \emph{intuitionistic} point of view on propositions and their proofs. There is nothing in this text that could support the claim that Kolmogorov's identifies here the  ``intuitionistic point of view'' with his own. The fact that in 1937 Kolmogorov did not convert himself into an intuitionist is independently supported by his 1936 Preface to Heyting's book \cite{Kolmogorov:1936} where Kolmogorov explicitly distances himself from the Mathematical Intuitionism (quoted in \textbf{2.5.} above); Kolmogorov's commentary of 1985  on his 1932 paper \cite{Kolmogorov:1991} quoted in  \textbf{2.1.}  is another evidence that his views on IL and on Mathematical Intuitionism did not change after 1932.

\subsection{Kolmogorov-Heyting Controversy}
Heyting's intuitionistic interpretation of the notion of proposition went through an interesting development before it achieved a stable form in his 1934 book \cite{Heyting:1931}, see also \cite{Sundholm:1983}. A part of this development involved an exchange of ideas between Heyting and Kolmogorov through their correspondence. Unfortunately, only three letters of Kolmogorov to Heyting first published by \cite[p. 14-16]{Troelstra:1990} (English translation \cite{Kolmogorov:1988}) are presently available. Heyting's replies to Kolmogorov are missing. Here we briefly overview this development and the letter exchange. As we shall see, Kolmogorov's view on IL equally went through a certain development before it stabilised in his 1932 paper \cite{Kolmogoroff:1932}.

In his \cite{Heyting:1930a} already quoted above Heyting introduces a distinction between proposition $p$ and proposition $+p$, which reads ``$p$ is provable''. He remarks that even if 

\begin{quote}
[f]ormulas $\vdash +p$ has exactly the same meaning as $\vdash p$ [\dots], [nevertheless] $p$ does not coincide with $+p$. \cite[p. 308]{Heyting:1998a}
\end{quote}
  
Heyting illustrates the distinction with the example of proposition "Every even number is a sum of two primes" (Strong Goldbach Conjecture or $SGC$ for short) and argues that judgements $\vdash \neg +SGC$ and $\vdash\neg\neg SGC$ are compatible. His argument goes as follows. $SGC$ is an expectation, taking an even number $k$ at random, to represent it as a sum of two primes. For every finite $k$ the possibility to do that is decided with a finite number of steps: one either obtains the wanted construction or shows that it is impossible.  $+SGC$, in its turn, ``requires a construction that gives us the decomposition [into two primes] for all even numbers at the same time''  \cite[p. 308]{Heyting:1998a} or, in other words, a general method $M$ that allows one to decompose into two primes any given even number.  
\footnote{Since $SGC$ was an open conjecture at the time of writing of Heyting's \cite{Heyting:1930a} and remains open today, the negation of $SGC$ is an example of what Brouwer calls a ``fleeing property'', which he defines as ``a property for which in the case of each natural number one can prove either that it exists or that it is absurd, while one cannot calculate a particular number that has the property, nor can one prove the absurdity of the property for all natural numbers'' \cite[p. 51]{Brouwer:1998}. Brouwer uses such \emph{fleeing properties} of natural numbers for constructing some of his \emph{weak counter-examples} to the unrestricted using LEM in mathematical proofs.}
Since the two expectations are different, it is possible that (i) the expectation for $+SGC$, i.e. the hypothesis that method $M$ exists, is lead to contradiction, and that (ii) the hypothesis that the decomposition fails at some $n$ is also lead to contradiction, which fulfils the expectation for $\neg\neg SGC$. As Heyting remarks in the conclusion of this argument   

\begin{quote}
[t]he difference between $p$ and $+p$ disappears if $p$ requires a construction.
\end{quote}

which shows that at this point he does not assume yet that \emph{every} proposition requires (or is an expectation of or intention toward) a construction (as he states later in 1934 \cite{Heyting:1934}, see below). 

A similar argument reappears in \cite{Heyting:1931}, see  \cite[p. 60]{Heyting:1931}. Here Heyting makes his first step toward \emph{abandoning} the $+$ operator. He first remarks that problems of form $+p$ always``require a construction '' (since $++p = +p$) and finally expresses doubts that the $+$ operator has any significance beyond a ``minimal practical'' one
\footnote{A similar point is made by Heyting in his letter to Oskar Becker of September 23, 1933, see \cite[p. 8]{Troelstra:1990} }. 
In the 1934 book \cite{Heyting:1934} Heyting abandons the  $+$ operator altogether and states that 

\begin{quote}
[e]ach mathematical proposition [$\dots$] is an intention towards a mathematical construction, which should satisfy certain conditions. \cite[p. 14]{Heyting:1934}, my translation from German
\end{quote}

By Sundholm's word, here Heyting ``finally commits himself'' \cite[p. 158]{Sundholm:1983}. As we shall now see Kolmogorov makes his commitment in \cite{Kolmogoroff:1932}, and this commitment is different.

The First letter of Kolmogorov to Heyting is dated October 12, 1931. The dates of the Second and the Third letters are missing but it is clear that the letters were sent in this chronological order (the First, the Second, an then the Third), that the Second letter has been sent before the beginning of the International Congress of Mathematicians in Zurich (September 12, 1932) and that the Third letter has been sent after Heyting's 1934 book appeared in press. This gives us a very approximate timeline for this correspondence.

In the \emph{First} letter Kolmogorov refers to Heyting's papers \cite{Heyting:1930a}, \cite{Heyting:1931} and objects to Heyting's aforementioned argument that $\vdash \neg +SGC$ and $ \neg\neg SGC$ are consistent: 

\begin{quote}
You consider [in \cite{Heyting:1930a}] as an example the statement ``Every even number is the sum of two primes''. At the same time, it is known that the formula $\vdash \neg\neg p \rightarrow p$ holds in this case both from the classical \emph{and from the intuitionistic} points of view. If it is established that $\vdash \neg\neg p \rightarrow p$, then we automatically have [in the given case] ''a construction that, for all even numbers, \emph{at once} gives us this decomposition.'' Therefore, $\vdash \neg\neg p \rightarrow +p$, and the case $\vdash \neg\neg p \wedge \neg +p$ is impossible. \cite[p. 89]{Kolmogorov:1988} translation corrected, my emphasis
\end{quote}

In modern terms Kolmogorov's objection can be formulated as follows. The statement of $SGC$ is expressed by a $\Pi^{0}_{1}$ formula (in the sense of Kleene-Mostowski arithmetic hierarchy), and for the formulas of this class the double negation elimination holds intuitionistically. I illustrate the proof with the example of $SGC$.  Since we have a general method to decide $G(n)$ for any $n$, formulas  $\forall n(G(n)\vee\neg G(n))$  and hence $\forall n(\neg\neg G(n)\rightarrow G(n))$ are justified intuitionistically. Since the universal quantifier distributes over the implication (both classically and intuitionistically), we get $\forall n(\neg\neg G(n))\rightarrow \forall n G(n)$. Finally, since $\forall x \neg\neg p(x) \rightarrow \neg\neg\forall x p(x)$ is an intuitionistic tautology we get  $\neg\neg\forall n G(n)\rightarrow \forall n G(n)$, i.e., $\neg\neg SGC \rightarrow SGC$. Notice that since the universal quantifier does not distribute over disjunction, formula  $SGC \vee \neg SGC$ cannot be proved by a similar argument. This latter formula is not intuitionistically valid because by the date of writing it remains unknown which of the the two disjuncts is true. The above argument is well-known and is not supposed to be original
\footnote{Notice that in the above quote (that dates to October 1931) Kolmogorov refers to the double negation elimination property of $SGC$ as a common knowledge. What is his possible reference?  A possible source is Hilbert's  Hamburg lecture of 1927 published in 1928 as \cite{Hilbert:1928}, English translation  \cite{Hilbert:1967}, where a distinction between ``real'' (i.e., variable-free and finitary decidable) and ``ideal'' (involving quantification of infinite domains) propositions is drawn, and a method of formal proof by ``adjoining'' ideal propositions to real ones (aka the ``method of ideal elements'')  is described. Hilbert illustrates this method with another $\Pi^{0}_{1}$ statement, namely, with the Last Fermat Theorem ($LFT$) by showing that given the consistency of formal arithmetic, the absurdity of  $\neg LFT$ is equivalent to $LFT$   \cite[p.471]{Hilbert:1967}. Even if this example is different (since $LFT$ is a negative proposition), a similar argument applies to all  $\Pi^{0}_{1}$ propositions including $SGC$ .}.   

Thus, anachronistically, it is clear that Heyting's argument was erroneous, and that Kolmogorov rightly pointed to Heyting's mistake. Historically, this episode is important because it could contribute to Heyting's decision to streamline his intended interpretation of propositions in IL by abandoning the $+$ operator (in  \cite{Heyting:1934}). But Kolmogorov's way out of the difficulty is different, as we shall now see. The First letter continues as this:

\begin{quote}
It seems to me that the problem here is not some defect of this particular example [of $SGC$]. Every ``proposition'' [$p$] in your [i.e., Heyting's] conception is, in my view, of one of the following two kinds:\\
($\alpha$) $p$ expresses the hope that, in some circumstances or other, some experiment will always give a definite result (for example, that the attempt to decompose any [given] even number $n$ into the sum of two primes gives a positive result [$\dots$]
\footnote{
As the example of $SGC$ example suggests, the condition that ``some experiment will always give a definite result'' should be read in the sense ``an experiment will give the \emph{same} expected result for all values of individual variables ranging over infinite domains''. 
}
. Naturally, every ``experiment'' must be realizable by means of a finite number of determined operations.\\
($\beta$) $p$ expresses the intention to find a certain construction.\\
We agreed that in case ($\beta$), the distinction between $p$ and $+p$ is inessential, but the proposition $\neg\neg p \rightarrow p$ must not be considered as
obvious. In the first case ($\alpha$), on the other hand, $p$ and $+p$ have different meanings, but we have $\vdash\neg\neg p \rightarrow p$ and $\vdash\neg\neg p \rightarrow +p$. [$\dots$] \\

I prefer to reserve the name of proposition [Aussage] only for propositions of the form ($\alpha$) and to call ``propositions'' of the form ($\beta$) simply problems [Aufgaben]. With the proposition $p$ are associated the problems $\neg p$ (to reduce $p$ to a contradiction) and $+p$ (to prove $p$). \cite[p. 89-90]{Kolmogorov:1988}, translation corrected
\end{quote}

As we see, at that point Kolmogorov still discusses problems and propositions in terms similar to Heyting's; in particular, in case ($\beta$) he accepts Heyting's identification of \emph{problem} $p$ with  \emph{proposition}  $+p$ ``$p$ is solvable''. But he also diverges here from Heyting by reserving the name of ``proposition'' only to case ($\alpha$). I disagree with van Atten who describes Kolmogorov's distinction between cases ($\alpha$) and ($\beta$) as a ``terminological matter'' and claims that this is what Kolmogorov meant when he made this distinction in his letter to Heyting \cite{vanAtten:2022}.  Kolmogorov's thinking behind this distinction, as I understand it, is this. 

Consider problem ``to decompose number 8 into two primes''. This is an ``intention to find a certain construction'', namely, to find a sum of two primes equal to 8. The problem is solved by pointing to the sum $5+3$. The $+$ operator in this case does not bring anything new: in order to show that the problem is solvable one only needs to present its solution. So in this case the distinction between proposition ``number 8 is decomposable into two primes'' and problem ``to decompose number 8 into two primes'' indeed appears as merely linguistic and logically sterile. However,  when one leaves the secure background of finitary reasoning,  the situation changes. Kolmogorov states (in case $\alpha$) that $SGC$ ``expresses the hope'' that the problem of decomposition into sum of two primes will be solved not only for number 8 but similarly for any other even positive integer. Notice that  Kolmogorov after Heyting explains here  $SGC$ in \emph{pseudo-constructive} terms, which involve infinite series of finitary constructive procedures 
\footnote{In his letter to Becker of September 23, 1933 Heyting describes such pseudo-constructive explanations of propositions ``non-constructive expectations'' (nicht-konstruktiven Erwartungen) \cite[p. 8]{Troelstra:1990}}.
Thus Kolmogorov's reason for distinguishing between cases  ($\alpha$) and ($\beta$), as I understand it, is to distinguish between properly constructive contents as in case ($\beta$) and pseudo-constructive ones as in case ($\alpha$). Another way out of the same difficulty could be, of course, to explain $SGC$ and other sentences in some constructively acceptable terms as in BHK-semantics. Heyting's refined interpretation of mathematical propositions in \cite{Heyting:1934} is a step in this direction. But this way out of the difficulty, as we shall shortly see, is not Kolmogorov's.

In his Second letter to Heyting (unfortunately, undated) Kolmogorov already judges about the same matter differently;
 
\begin{quote}
I have thought about your example of sentence \\``For all $i$ we have $a_{i} < b_{i}$''.\\  In general,  let $x$ be a variable and $P(x)$ a problem depending on this variable. The ``hope'' of finding, for any $x$, a solution of $P(x)$  is, in my terminology, neither a problem nor a proposition
\footnote{In the original German the sentence reads: Die ``Hoffnung'' f\"ur jedes $x$ eine L\"osung der Aufgabe $P(x)$ zu finden ist in meine Terminologie weder eine ``Aufgabe'' noch eine ``Aussage''. Once again Kolmogorov's wording is a bit confusing here. As it stands the argument may appear unsound. The ``hope'' of solving the problem ``prove that $a_{i} < b_{i}$'' is unproblematic for any given $i$ (assuming that series of numbers $\{a_{i}\}$ and $\{b_{i}\}$ are well defined). What is problematic is to verify this property for \emph{all} values of $i$ as in Heyting's example. So I read the cluase ``finding, for any $x$, a solution of $P(x)$'' as ``finding a solution of problem $P(x)$ for all $x$''.  
}. 
It would be very interesting to learn whether you associate with this hope a positive expectation that, for any $x$, $P(x)$ will actually be solved (when and by whom?). If this expectation is not implied, then I think we are approaching a naive non-intuitionistic understanding of the assertion ``$P(x)$ is soluble for all $x$''. \cite[p. 90]{Kolmogorov:1988}, translation corrected
\end{quote}

When Kolmogorov says that pseudo-constructive explanations suggested earlier by Heyting in \cite{Heyting:1930a} are not appropriate either to problems or to propositions, he implies that they are ill-formed. Kolmogorov's rhetorical questions (``solved when and by whom?'') are supposed to substantiate the claim by stressing that the ``hope'' to solve problem of form $P(x)$ for \emph{all} values of $x$, where $x$ ranges over an infinite domain, cannot be fulfilled in principle. Here Kolmogorov sounds as an ultra-intuitionist who does not want to accept the mathematical infinity in any form including its constructive form of  \emph{potential} infinity. But this is certainly not his position, as one can judge on the basis of the above analysis of Kolmogorov's 1932 paper \cite{Kolmogoroff:1932}, and also on the basis of his wider mathematical practice not constrained by the limits of constructive mathematics. Kolmogorov's way out of the difficulty is to reserve the intuitionsitic style of mathematical reasoning (including IL as its formalisation) to problem-solving and at the same time allow classical reasoning in proofs of propositions (theorems). This allows Kolmogorov to avoid pseudo-constructive reasoning both in explaining propositions and in explaining problems. In this setting \emph{proposition} $SGC$ admits a classical explanation in terms of (the infinite set of) ``all'' even natural numbers, while the associated \emph{problem}  ``to prove $SGC$'' admits a \underline{constructive} explanation, which doesn't relay on the doubtful pseudo-constructive ``hope'' that the natural numbers can be \emph{all} somehow checked one by one but relies instead on what can be \underline{classically} admitted as a proof of $SGC$. How exactly this interplay between the constructive and the classical mathematical reasoning works Kolmogorov leaves in \cite{Kolmogoroff:1932} unexplained; he only insists here that the intuitionistic way (and form) of mathematical reasoning applies only to problems (including problems of form ``to prove proposition $p$'') but not to propositions. As we have shown in \textbf{2.3.}, in the domain of problems Kolmogorov applies the intuitionistic requirements uncompromisingly asking the reader to \emph{have really solved} a given problem before moving on to the next problem \textemdash\ rather than just assume the given problem to be solvable or be earlier solved by other people. Melikhov's work \cite{Melikhov:2022a}, \cite{Melikhov:2022b} is an attempt to formalise the interplay between classical and intuitionistic reasoning using Kolmogorov's insights.              

Heyting also gets rid of pseudo-constructive explanations but he does it in a different way, which, unlike Kolmogorov's way, can be called properly intuitionistic. He pursues his strategy to explain mathematical propositions in terms of problems, expectations (``hopes'') and intentions and abandons the $+$ operator along with pseudo-constructive explanations, which motivate the introduction of this operator. So he comes in \cite{Heyting:1934} to his mature (``committed'' by Sundholm word) interpretation of mathematical propositions that leaves no room for distinguishing between problems and propositions in a non-trivial way. According this new view any admissible proof of a proposition is a genuine construction \textemdash\ rather than a pseudo-construction like an accomplished check of every even number showing that every even number is decomposable into a sum of two primes. What qualifies here as a construction in the general case remains to be further specified but in case of $SGC$ an obvious candidate is a general method or algorithm that inputs an even number and outputs its decomposition into two primes (but not just checks if the given even number is so decomposable) and that somehow makes it clear that it cannot fail whatever even number is input. 

Since by the time of writing no proof of $SGC$ is known it remains today an open problem. Both Kolmogorov and Heyting would readily agree with that. But they would interpret this sentence differently; and their different interpretations would imply different views on what could count as a solution of the problem. For Kolmogorov the open problem is ``to prove $SGC$'', in symbols $+SGC$ but not $SGC$ as such (as a conjectured proposition). For Heyting (as Kolmogorov perfectly explains in his note of 1937 \cite{Kolmogoroff:1937}) $SGC$ in its present conjectural status is itself  an open problem (without the + operator). A major consequence of this difference in views is that for Heyting only a constructive proof may qualify as a solution of $SGC$ while Kolmogorov is in a position to admit as the wanted solution a non-constructive proof of this conjecture as well 
\footnote{A modal interpretation of IL proposed by Kurt G\"odel in 1933 \cite{Godel:1933a} (English translation \cite[301-303]{Godel:1986})  has a similar feature. Here propositional variables of IL are interpreted as modalised \emph{classical} propositions of form $B(p)$ read ``proposition $p$ is provable'' (Compare with Heyting's + operator.) G\"odel shows that, on the pain of contradiction, the provability modality $B$  should be understood in the ``absolute'' sense rather than in a special sense as the provability in some fixed formal system. Thus in G\"odel's 1930 interpretation of IL classical proofs  fall under the relevant concept of provability. Like Kolmogorov's interpretation G\"odel's interpretation of IL does \emph{not} qualify as intuitionistic (if one takes \cite{Heyting:1934} as the standard of being intuitionistic). 
},
\footnote{The \emph{Third} letter of Kolmogorov to Heyting \cite[p. 91]{Kolmogorov:1988}, also undated, contains no material relevant to the present discussion. In this letter Kolmogorov thanks Heyting for sending him his book of 1934 \cite{Heyting:1934} and points to Heyting that Kurt G\'odel's results published in his 1933  paper on the intuitionsitc arithmetic \cite{Godel:1933} (English translation \cite[287-295]{Godel:1986}) are ``very close'' to Kolmogorov's results published in his Russian paper \cite{Kolmogorov:1925} back in 1925 (see \textbf{2.1} above). }.

.

\section{The Legacy of Kolmogorov's Calculus of Problems beyond the BHK-interpretation} 
The aim of this Section is twofold. First, we provide a brief overview of some later works motivated by Kolmogorov's CP as described in his 1932 paper  \cite{Kolmogoroff:1932} independently of the related Heyting's works and hence of the BHK-interpretation of IL. Second, we explore a new interesting perspective on the Kolmogorov-Heyting controversy concerning the distinction between problems and propositions (as described in the last Section) provided by Homotopy Type theory and Univalent Foundations of mathematics.    
    
\subsection{Works motivated by Kolmogorov's Calculus of Problems}
Kolmogorov's 1932 paper \cite{Kolmogoroff:1932} has been published in German in an influential mathematical journal of the time and immediately reached its intended audience which included Arend Heyting, Kurt G\"odel (who refers to this Kolmogorov's paper in his \cite{Godel:1933a}) and many other researchers in the field of logic and foundations of mathematics. So it is hardly possible to trace here all later developments where Kolmogorov's logical ideas played a role; moreover so, after  these ideas were combined in 1980s by Troelstra  with Brouwer's and Heyting's ideas within the so-called BHK-interpretation of intuitionistic logic \cite{Troelstra:1990}, which became very popular. Our task would become even more difficult if we try to take into account results and developments, which are strongly related to Kolmogorov's interpretation of IL theoretically without being motivated by Kolmogorov's work in the real history. The \emph{realisability interpretation} of intuitionistic arithmetic due to Stephen C. Kleene is a case in point  \cite{Kleene:1945}
\footnote{
In his autobiographical paper \cite{Kleene:1973} Kleene claims that neither Heyting's 1934 proof interpretation of IL, nor Kolmogorov's 1932 interpretation (about which Kleene was aware via the same source \cite{Heyting:1934}), played a role in developing his original realizability interpretation of inuitionistic sentences. About Kolmogorov's interpretation he says that it ``failed to help me in any way of which I'm conscious'', see \cite[p. 100, footnote 6]{Kleene:1973}. 
}.
For the above reasons we limit the following review only to works influenced by Kolmogorov's 1932 paper explicitly and directly, including the works of some Kolmogorov's students. 

The chronologically first systematic attempt to develop Kolmogorov's calculus has been made by Kolmogorov's Ph.D. student Yuri Tikhonovich Medvedev in his dissertation ``On the degrees of difficulty of mass problems'' defended in 1955, see \cite{Medvedev:1955} for a summary and \cite{Mostowski:1956} for its review in English. The key idea here is to represent a given problem $P$ by set $S$ of arithmetical functions $f : \mathbb{N}\rightarrow \mathbb{N}$, which can be thought of as potential solutions of $P$. Medvedev calls \emph{mass problems} those problems, which are so representable. $P$ is algorithmically solvable when $S$ contains a \emph{recursive} function. Further, Medvedev defines (using recursive functionals) a relation of reducibility between mass problems, constructs on this basis lattice $\Omega$ of their ``degrees of difficutly'' (so that all pairs of mutually reducible problems have the same degree), and, finally, shows that $\Omega$ models IL. In 1956 Albert Ivanovitch Muchnik showed that Medvedev's ``degrees of difficulty'' were, generally, not algorithmically comparable \cite{Muchnik:1956} and later proposed a weaker non-constructive reduction relation between Medvedev's mass problems that didn't require exhibiting an explicit recursive functional \cite{Muchnik:1963}, English translation \cite{Muchnik:2016}. In the late 2000s these works by Medvedev and Muchnik anew attracted a considerable attention \cite{Simpson:2008} (see also reference therein), and in 2016 a modern computational logical framework based on Medvedev's and Muchnik's ideas was proposed by S.S. Basu and S.G. Simpson under the name (and mathematical form) of \emph{Muchnik topos} \cite{Basu&Simpson:2016}. 

In 1960s Medvedev published a series of papers  \cite{Medvedev:1962}, \cite{Medvedev:1963}, \cite{Medvedev:1966} where he explored a somewhat different approach to formalising Kolmogorov's concept of problem.  Here Medvedev focuses on problems that admit only a finite number of possible solutions without assuming that these possible solutions are arithmetical functions. Such problems Medvedev calls \emph{finitary}. He shows, in particular, that IL is incomplete with respect to this new semantics and proposes some syntactic modifications of IL aiming to developing a new logical calculus. This line of research was further pursued in 1970s and 1980s by Medvedev's colleague Dmitri Pavlovitch Skvortsov who proposed some generalisations of Medvedev's approach including a generalisation to the infinitary case \cite{Skvortsov:1979}, \cite{Skvortsov:1979a}. A concise but more detailed review of these works is found in \cite{Uspenskii&Plisko:1991}. Commenting on Medvedev's and Skvortsov's works referred to above Uspenski and Plisko remark that 

\begin{quote}          
The new interpretation of intuitionistic logic [due to Kolmogorov], free of philosophical concepts of intuitionism, made it meaningful to investigate the logic as a calculus of
problems.  \cite[p. 460]{Uspenskii&Plisko:1991}
\end{quote}

Indeed Medvedev, Muchnik, Skvortsov and many other Russian mathematicians who developed Kolmogorov's logical ideas were hardly interested in philosophical foundations of logic and mathematics but contributed instead to the Recursion theory, theory of Algorithms, and, more generally theory of Computation. Their work should be understood in the context of contemporary developments in these and other related areas of mathematics and computer science, about which the Russian researches were usually aware albeit often with certain delays. The aforementioned result of A.A. Muchnik concerning the Medvedev lattice is remarkable in this respect. During the same year of 1956 the same result was obtained by American mathematician and theoretical physicist Richard M. Friedberg in a different setting, namely,  as a response to a problem in the Recursion theory formulated by Emil Post in 1944 \cite{Friedberg:1956}. The result is known today under the name of Friedberg-Muchnik theorem  \cite[p.253]{Kozen:2006}.  Muchnik was aware about Friedberg's publication by 1959 \cite{Muchnik:1959}; I could not find it out when Muchnik's Russian publication became first known to Friedberg. This example demonstrates how Kolmogorov's interpretation of IL helped to establish theoretical connections between the Intuitionistic Logic and other related areas of the contemporary mathematical research, which eventually gave birth to Computer Science. 

Kolmogorov's interpretation of IL also plays a major role in the work of his former student Sergei Nikolaevitch Artemov. Since late 1980s Artemov is developing a research program in Epistemic Logic that started from studying various Logics of Provability \cite{Artemov&Beklemeshev:2004} and eventually developed into the original concept of Justification Logic \cite{Artemov:2008}, \cite{Artemov&Fitting:2019}. A starting point of this program was the idea to integrate concrete individual proofs (aka proof terms) into the G\"odel-style modal Provability Logic  \cite{Artemov:1994a}. Describing his motivations Artemov eventually refers to the BHK-semantics but he makes it clear that his approach is based on a combination of Kolmogorov's interpretation with G\"odel's 1933 provability interpretation of IL \cite{Godel:1933a} (that Artemov calls the ``Kolmogorov-G\"odel'' approach \cite{Artemov:2004}) rather than on a combination of Kolmogorov's interpretation with Heyting's, which underpins the notion of BHK semantics.  Kolmogorov's Calculus of Problems equally plays a major motivating role in the research of Artemov's former student Giorgi Japaridze working on a general Computability Logic \cite{Japaridze:2003}, \cite{Japaridze:2003}, which he sees as ``a justification and hence a materialization of Kolmogorov’s known thesis '[IL] = logic of problems' '' \cite[p. 77, abstract]{Japaridze:2007}. Artemov's and Japaridze's research demonstrate, once again, the fact that the notion of BHK interpretation, however useful for many theoretical purposes, should not be seen as a vehicle that fully and faithfully embodies the ideas of Kolmogorov's 1932 logical paper \cite{Kolmogoroff:1932}. 

We conclude this short list of works motivated by Kolmogorov's CP with the recent work of Sergei Melikhov who proposed a \emph{Combined Logic of Problems and Propositions} (QHC) intended as a realisation of Kolmogorov's idea of  ``unified logical apparatus dealing with objects of two types \textemdash\ propositions and problems. \cite[p.452]{Kolmogorov:1991}'', \cite{Melikhov:2022a}, \cite{Melikhov:2022b}, \cite{Melikhov:2023}. Here is a very informal description of QHC. It comprises (i) a copy of classical predicate calculus (QC) interpreted as usual, (ii) a copy of intuitionistic predicate calculus (QH) interpreted after Kolmogorov in terms of problems and their solution, and (iii) two operators denoted ? and ! with the following intended meaning:

\begin{itemize}
\item given \emph{proposition} $p$, formula $! p$ denotes the  \emph{problem} ``to prove $p$'';
\item given \emph{problem} $\alpha$, formula $? \alpha$ denotes the  \emph{proposition} ``$\alpha$ has a solution''. 
\end{itemize}

Notice that $?! p$ says ``the problem of proving proposition $p$ has a solution'', i.e., ``$p$ is provable'', and $!? \alpha$ says ``to prove that problem $\alpha$ has a solution'', i.e., ``to prove that $\alpha$ is solvable''. 

In addition to the standard axioms of QC and QH the combined calculus QHC comprises five axioms: 
\begin{itemize}
\item $?! p \rightarrow p$ (``$p$ is provable'' implies $p$);
\item $\alpha \rightarrow !? \alpha$ (solution of problem $\alpha$ solves the problem ``to prove that  $\alpha$'' is solvable'');
\item $! (p \rightarrow q) \rightarrow (! p \rightarrow ! q)$ (a reduction of proof of $q$ to proof of $ p$ reduces to proving implication $(p \rightarrow q)$ );
\footnote{Beware of an abuse of notation in this and the following axiom: we use here the same symbol for the classical implication and for the intuitionistic implication (interpreted as reduction of problems).}
\item $? (\alpha \rightarrow \beta) \rightarrow (? \alpha \rightarrow ? \beta)$ (if the problem of  reducing $\beta$ to $\alpha$ is solvable then the solvability of $\alpha$ implies the solvability of $\beta$);
\item $\neg ! \bot$ (the falsity has no proof);
\end{itemize} 

and two additional rules of inference applying across QC and QH:

\begin{prooftree}
\hypo{p} 
\infer1[(CH)]{! p}
\end{prooftree}

and

\begin{prooftree}
\hypo{\alpha} 
\infer1[(HC)]{? \alpha}
\end{prooftree}
 
 \bigskip
Since these rules connect the classical and the intuitionistic fragments of  QHC, there is a choice between classical and intuitionistic interpretations of these rules themselves,  which can be understood as implications in the corresponding \emph{meta}logic). Melikhov assumes that the metalogic of QHC is intuitionistic and interprets both rules using Kolmogorov's interpretation of the intuitionistic implication. Thus rule (CH) reads (informally):

\begin{quote}
 There is a method to find, for each formula F, a method deriving from a proof of the assertion that each of the propositions instantiating F is true a solution of the problem of proving (by a general method) all propositions instantiating F. 
\end{quote}

 while rule (HC) reads:
 
\begin{quote}
There is a method to find, for each formula $\Phi$, a method deriving from any solution (by a general method) of all
problems instantiating $\Phi$ a proof of the assertion that each of the problems instantiating $\Phi$ has a solution.
\footnote{
Cf. \cite[section 3.2.]{Melikhov:2022b}. The above interpretations of rules (CH) and (HC) are communicated by Sergei Melikhov to the author privately. 
}
\end{quote}  
 
QHC admits both problem-based and proof-based interpretation and has interesting mathematical properties.  The most remarkable property is the Galois connection between the Lindenbaum posets of equivalence classes QH-formuals and QC-formulas provided by ? and ! operators: 
\begin{center}
$? \alpha \rightarrow p $ if and only if $\alpha \rightarrow ! p$
 \end{center}
 
Extending QC with the modal operator $\square := ?!$ returns the familiar modal logic S4 while extending QH with the dual operator $\nabla := !?$ brings a version of the Intuitionistic Epistemic Logic studied by Artemov and Protopopescu \cite{Artemov&Protopopescu:2014} (as well as the Lax Logic of Fairtlough and Walton \cite{Fairtlough&Walton:1997} and the Russell-Prawitz modal logic of Aczel \cite{Aczel:2001}, see \cite[p. 22]{Melikhov:2022a} for further details). QHC also admits a topological semantics  \cite[section 4]{Melikhov:2022b} as well as a Kripke-style semantics described by Anastasia Onoprienko, which she used for establishing the completeness of QHC with respect to this semantics \cite{Onoprienko:2020},\cite{Onoprienko:2022}. 

QHC along with its intended interpretation is fully compatible with Kolmogorov's interpretation of IL presented in \cite{Kolmogoroff:1932} and with Kolmogorov's philosophical justification of this interpretation analysed in the present paper. In that sense QHC certainly qualifies as a realisation of Kolmogorov's project. However, QHC obviously goes beyond the direct indications found in Kolmogorov's writings. In particular, in Kolmogorov's writings I cannot see any pointer to Melikhov's ? operator that associates with a given problem $\alpha$ a proposition saying that $\alpha$ is solvable. While the idea to associate a problem with a given proposition (Melikhov's operator !) is immediately read off from Kolmogorov's  examples, this is not the case of the ? operator.  This remark suggests that Kolmogorov's ideas found in his 1932 paper and in related writings may also admit different interpretations and different developments. 

\subsection{Calculus of Problems and Univalent Mathematics} 
Here I show how the Homotopy Type theory (HoTT) and Univalent Foundations of Mathematics (UF) help to justify Kolmogorov's distinction between problems and propositions  in a purely constructive setting without using the classical Frege-style notion of proposition that has a truth-value independently of one's knowledge of this proposition (and hence of one's capacity to prove or disprove it). Following Heyting, we assume here the intuitionistic notion of proposition as ``an intention towards a mathematical construction, which should satisfy certain conditions'' \cite[p. 14]{Heyting:1934}. But we do this not without a further ado. As we shall shortly see, in HoTT, the notion of mathematical construction (which satisfies certain conditions) admits further qualifications. In what follows I explain (or rather motivate) this point very informally without referring to HoTT and then provide some mathematical details.

Consider the classical geometrical problem (Proposition 1.1 of Euclid's \emph[p. 8]{Elements}):

(E) : \emph{to construct a regular triangle on a given side by ruler and compass}

and the corresponding proposition

(H) : \emph{for any given straight segment, there exists a regular triangle constructible on this segment (taken as its side) with ruler and compass}

Let us assume that (H) is interpreted accordingly to Heyting as an ``intention towards a construction'', namely a construction of regular triangle by ruler and compass satisfying the aforementioned condition (that one side of the constructed triangle should coincide with the given line).  According to Heyting (E) and (H) are just two different wordings of the same proposition. However, without leaving the constructive ground, one can  interpret (E) in a stronger sense than (H) assuming that (E) amounts to construction of a full-fledged geometrical object \underline{provided with its identity conditions}\footnote{Recall of Quine's motto ``no entity without identity''.} while (H) does not involve this latter requirement. The talk of identity conditions in case (E) implies a need to consider the \emph{type} of regular triangles and then specify how individual \emph{tokens} (aka terms) of this type, i.e., the individual triangles, are mutually identified and distinguished. We shall see that in HoTT such identity conditions are specific to types and, generally, very non-trivial. As it is usual in the constructive mathematics (and in accordance to Heyting's interpretation of propositions) one assumes here that (H) can be satisfied with any schematic construction, which shows that an object with the wanted properties is constructible \textemdash\ without such additional identity considerations. The distinction between problems and propositions is now the following: while a problem requires establishing specific identity conditions of a constructed objects, a proposition requires only a schematic construction (assuming that the relevant identity conditions are trivial).  

Such a specific reading of (E) may seem artificial but the fact that Euclid introduces the ``given side'' in his formulation of the problem (rather than simply asks to produce a regular triangle from the scratch using some specified constructive procedures) suggests that he indeed had some identity considerations about geometrical objects in his mind. Remark that Euclid's construction of regular triangle by ruler and compass produces congruent triangles if and only if the given straight segments are congruent, and that the congruence of straight segments is easier to check than that of triangles. So the given straight segment determines the ``individuality'' of the produced triangle up to congruence. \footnote{Notice that unless one breaks the symmetry by making an arbitrary choice at an earlier stage of Euclid's construction (as Euclid himself does), the construction of Proposition 1.1 produces not one triangle but two congruent triangles.}. It is also suggestive to think of ``giving'' a straight segment as an individual cognitive act that determines the identity conditions of the following construction that results in the wanted triangle, which, in its turn, is followed by the final verification that the constructed object is as required. 

At the same time it is clear that in Euclid's \emph{Elements} the identity conditions of geometrical objects are not rigorously fixed. In some contexts congruent figures can be thought of as the ``same'' but in some other contexts they should be thought of as different (for example, the three sides of a regular triangle are congruent but not the same)
\footnote{Euclid's notion of \emph{equality}, which is implicitly defined by his Common Notions aka Axioms, does not solve this problem. It is not what a modern reader normally expects: in case of two plane figures $F, G$ Euclid's equality translates into modern terms by saying that  $F$ and $G$ have equal areas. Euclid's Axiom 4 says that congruent objects are equal (in the above sense) but the converse does't hold. Thus Euclid's mathematics provides no definite answer to the question of whether two given geometrical objects are the same or not. On may even argue that this question is simply irrelevant in this theory (unlike, say, the question of whether or not two given figures have the same area).}.  
It is more surprising that the same traditional looseness about the identities of mathematical objects is also characteristic for the 20th century mathematics. This is in spite of the fact that     
in the modern set-theoretic setting the identity problem was supposed to be fixed once and for all with the formalised identity relation making part of the underlying logical machinery (the classical predicate calculus with identity in case of ZF). Here a single identity relation is supposed to serve universally for all mathematical needs. As  remarks Frege 
\begin{quote}
Identity is a relation given to us in such a specific form that it is inconceivable that various forms of it should occur'' \cite[p. 254]{Frege:1903}, my translation from German
\end{quote}
and ZF as well as many other formal theories built during the 20th century straightforwardly implement this Frege's idea. 

As it is well-known today, this universal solution of the identity problem turned out to be very unsatisfactory in the mathematical practice and led, in particular, to the so-called \emph{Benacerraf Identification Problem} widely discussed by philosophers of mathematics during the last decades \cite{Benacerraf:1965}. The problem in its general form is the following. ZF justifies the principle of \emph{Indiscernibility of Identicals} (sometimes called the \emph{Leibniz Law}) according to which identical objects have identical properties, in symbols

\begin{equation}\label{eqn:InId}\tag{InId}
(x=y)  \rightarrow \forall P  \ldotp (P(x) \leftrightarrow P(y))
\end{equation}

Here the Indiscernibility of Identicals is formulated as an axiom of Second Order logic but it can be also formulated as an axiom schema in the First Order logic. In the ZF-based mathematics every informal mathematical property $P$ of certain objects $x$ of type $X$ is represented via the Separation Scheme of ZF as a subset $P\subseteq X$, so that for all $x \in X$, $x\in P$ just in case $P(x)$ holds, that is, $x$ has property $P$ \footnote{We deliberately abuse here the notion by denoting informal mathematical objects and their properties, on the one hand,  and their corresponding set-theoretic representations, on the other hand, by the same symbols.} In this more specific setting the  Indiscernibility of Identicals is expressed by formula

\begin{equation}\label{eqn:InIdZF}\tag{InIdZF} 
(x=y)  \rightarrow \forall p \ldotp (x\in p \leftrightarrow y\in p)
\end{equation} 

which follows immediately from the usual axioms of identity making part of ZF's logical machinery. In versions of ZF without identity the right part of (InIdZF) is used as a definiens for \emph{defining} the identity relation $x=y$.  Notice that in this discussion we wholly leave aside the converse principle of \emph{Identity of Indiscernibles}, which is more controversial.  

(InId) is a basic logical principle without which the concept of identity hardly makes sense. The problem is that the identity relation = of ZF turns out to be wholly irrelevant in mathematics outside ZF itself. Say, for group theorists the ``right'' notion of identity of algebraic groups is the  \emph{isomorphism} of groups, so they want the  Indiscernibility of Identicals principle in the following form, which Ahrens and North call an \emph{equivalence principle} \cite{Ahrens&North:2019}): 

\begin{equation}\label{eqn:InIdGr}\tag{InIdGr}  
 (G \cong H)  \rightarrow \forall P \ldotp(P(G) \leftrightarrow P(H))
\end{equation} 

where $G, H$ are isomorphic groups and $P$ ranges over all \emph{group-theoretic} properties. But there is no easy way \textemdash\ and apparently no way at all  \textemdash\ to upgrade or specify (InIdZF) to (InIdGr) because it is not known in advance which properties expressible in the language of ZF qualify as group-theoretic and which don't! This is why in Group theory the identity relation of ZF turns out to be simply irrelevant, and group theorists find themselves in the same situation as Euclid as far the issue of identity of mathematical objects is concerned. The set-theoretic language helps group theorists to formulate the basic definition of group and prove many theorems but it is not helpful (if not misleading) for tackling the identity issue.

\emph{Mutatis mutandis} the same can be said about the Ring theory or any other area of the ``structuralist'' mathematics where the relation of isomorphism between mathematical structures has the same or similar role. Benacerraf points to Zermelo's and von Neumann's ordinals as two equally good candidates for representing natural numbers stressing the fact that their shared structure cannot be properly identified by the standard means \cite{Benacerraf:1965}. Notice that the the wanted ``structural'' identity concept is type-dependent. Indeed, the  Indiscernibility of Identicals principle for rings

\begin{equation}\label{eqn:InIdR}\tag{InIdR}  
 (R_{1} \cong R_{2})  \rightarrow \forall P \ldotp(P(R_{1}) \leftrightarrow P(R_{2}))
\end{equation}    

has the same form as (InIdGr) but the two isomorphisms relations are not the same in both cases, and the corresponding classes of concerned properties also differ. Moreover, the isomorphism of set-based mathematical structures is not always the ``right'' notion of identity. There is a general consensus that in Category theory the right notion of identity is the \emph{category-theoretic equivalence}, which is weaker than isomorphism of categories. The rise of category-theoretic mathematics during the second part of the 20th century further highlighted the identity problem in mathematics but once again didn't provide a satisfactory general solution \cite[ch. 6]{Rodin:2014}.

The concept of Univalent Foundations of mathematics (UF) was first officially presented by Vladimir Voevodsky in 2010 in his lectures \cite{Voevodsky:2010},\cite{Voevodsky:2010a},  see also his earlier 2006 lecture \cite{Voevodsky:2006} where the same idea is introduced  under a different name. The theoretical basis of UF is HoTT, which is an interpretation of Intuitionistic Type theory due to Per Martin-L\"of (MLTT) in terms of Homotopy theory \cite{Grayson:2018}. Voevodsky's \emph{Univalence Axiom} (UA) implements in HoTT a very general form of \emph{equivalence principle}, that is, a far-reaching generalisation of (InId) \cite{Ahrens&North:2019}. 
On this basis UF provides a rigorous formal account of identity of mathematical objects, which is very unlike the standard account provided by ZF and other first-order theories. In the UF the structure of identity of given mathematical object $x$ depends on the corresponding type $X$ of objects and in a sense makes part of the construction of $x$ itself. This may appear less surprising when one takes into account the fact that in order to construct a mathematical structure $x$ of type $X$ \emph{up to isomorphism} the appropriate isomorphisms of structures of type $X$ also need to be constructed. So the idea that the identities of mathematical objects are constructed along with these objects themselves (as this happens in case of social groups and some other artefacts) is not, after all, so counter-intuitive as it may first seem
\footnote{
Let me also illustrate this point with the above elementary example of regular triangle. Suppose we want to construct such a triangle \emph{up to congruence}. For our purpose it is essential to think of congruence not only as a relation applicable across the Euclidean plane but also as a concrete invertible map $c : x \xrightarrow[\sim]{} y $ that maps given figure $x$ into a congruent figure $y$. The standard name for such a map is \emph{isometry} (distance)preserving map). The \emph{relation} of congruence $x \cong y$ holds when there is an isometry between figures $x$ and $y$. Now observe that there are 3! = 6 different isometries (including the trivial one) mapping regular triangle $ABC$ onto itself, which form symmetric group $S_{3}$ (the group of permutations of letters $A, B, C$). Three of these six isometries preserve orientation of the plane and form cyclic subgroup $C_{3} \subset S_{3}$. If we now take a generic scalene triangle $DEF$ we can see that there is a single isometry of  $DEF$ onto itself, namely the trivial one. Now we can see that the conventional expression ``the same up to congruence'' can (and arguably should) be understood  in case of $ABC$ and in case of $DEF$ in different ways \textemdash\ because these triangles are of different types, in which isometries apply in different ways. The difference disappears only when the group structure of isometries is ignored, and congruence is treated in the usual way as a mere relation. But since one assumes that congruence serves as the identity condition of a geometrical object, it is not unreasonable to further assume that the corresponding groups of isometries constitute a finer structure of the same identity concept. Since we have here a setting where a regular triangle is \emph{constructed} (rather than simply given) a similar constructive mode of reasoning may also extend to its groups of symmetries $C_{3} \subset S_{3}$ (albeit the ruler and compass are not quite appropriate instruments for building these groups). }. 

As its very name clearly indicates, the \emph{Intuitionistic Type Theory} was conceived of by Per Martin-L\"of as a novel formal carrier of the intuitionistic mathematical reasoning.  This calculus applies the idea of Carry-Howard Correspondence and is particularly apt for computational implementations. The first version of Intuitionistic Type theory was presented by Martin-L\"of in autumn of 1970 in a seminar lecture. This first version of MLTT was significantly revised after a critique of Jean-Yves Girard who showed its inconsistency. In what follows we refer to the 1980 corrected version of MLTT published as \cite{Martin-Lof:1984}, which remained by and large stable ever since. For a detailed history of MLTT and its philosophical motivation and background see  \cite{Sundholm:2012}. 

A key logical concept of MLTT is that of \emph{judgement}. A basic form of judgements (one of four forms listed in \cite[p. 5]{Martin-Lof:1984}) is $x : X$ which reads ``$x$ is a term of type $X$''. According to Martin-L\"of formula $x : X$  admits the following contentual interpretations (or, more precisely, ``explanations'') \cite[p. 5]{Martin-Lof:1984}: 

\begin{enumerate}
\item $x$ is an element of set $X$
\item $x$ is a proof (witness, evidence) of proposition $X$ 
\item $x$ is a method of fulfilling (realising) the intention (expectation) $X$
\item $x$ is a method of solving the problem (doing the task) $X$ 
\end{enumerate}

Notice that the combination of interpretations (2) and (3) is close to Heyting's interpretation of mathematical propositions as `an intention towards a mathematical construction, which should satisfy certain conditions'' \cite[p. 14]{Heyting:1934}; a difference being that Martin-L\"of  distinguishes between propositions and judgements \cite[p. xx-xxi]{Sundholm:2012} and avoids using the term ``construction'' (apparently, only for the reason of parsimony). Interpretation (4) clearly points to Kolmogorov's problem interpretation of IL as presented in his \cite{Kolmogoroff:1932}. Applying Kolmogorov's interpretation along with Heyting's Martin-L\"of thinks of the two interpretations as complimentary (in line of the BHK-interpretation, which has been first formulated about the same time): the idea here is that (1)-(4) explain in different words one and the same fundamental logical concept. Thus vis-\`a-vis the Kolmogorov-Heyting controversy as described in the present paper (in \textbf{3.1.} above), Martin-L\"of is on Heyting's side: proposition $X$ with its proof $x$ and the corresponding problem $X$ with its solution $x$ are treated by Martin-L\"of as the same judgement $x : X$; Kolmogorov's distinction between problems and propositions is not supported by MLTT at the formal level.       

Interpretation (1) is Martin-L\"of's original; it connects the BHK semantics to Set theory (and set-based mathematics) and to Type theory. This is how Martin-L\"of explains the equivalence of his interpretations (1) and (2):

\begin{quote}
If we take seriously the idea that a proposition is defined by lying down how its canonical proofs are formed [$\dots$] and accept that a set is defined by prescribing how its canonical elements are formed, then it is clear that it would only lead to an unnecessary duplication  to keep the notions of proposition and set [$\dots$]  apart. Instead we simply identify them, that is, treat them as one and the same notion.  \cite[p. 13]{Martin-Lof:1984}
\end{quote}

The equivalence of (1) and (2) allows one to identify a proposition with a set of its proofs. Martin-L\"of doesn't comment separately on the equivalence between (1) and (4) but since all interpretations (1)-(4) are supposed to be in the same sense equivalent we are also in a position to identify a problem with a set of its solutions. Recall from \textbf{4.1.} that such a set-theoretic representation of problems was independently used by Yu. T. Medvedev in 1960s.

As we shall now see, the homotopy-theoretic interpretation of MLTT (that is, HoTT) imposed a significant modification of its intended interpretation outlined above. Unlike the intended interpretation, the modified interpretation supports Kolmogorov's idea according to which propositions and general problems should be formally distinguished. It is remarkable that this modification is made for a purely mathematical reason but not as an attempt to implement this or that philosophical idea about logic and mathematics. 

Let $x=_{X}y$ be a type that in line of (2) is interpreted as a proposition that says that two terms $x,y : X$ of type $X$ are equal (i.e., are the same). Let $x_{1},y_{1} : x=_{X}y$ are two terms of this identity types that according to (2) we interpret as two proofs (witnesses) that the proposition is true. Under the homotopical interpretation the underlying type $X$ (as well as any other type) is interpreted as a \emph{space} in a sense sufficient to support the homotopy theory,  which is the case of a topological space. Terms $x,y$ are interpreted as \emph{points} of this space, and the identity proofs $x_{1},y_{1}$ are interpreted as continuous \emph{paths} between points $x,y$. (If there is a continuous path between points $x,y$ then this path can be ``shrunk'' or ``contracted'' showing that  the two points are indeed the same.)

The syntax of MLTT allows one to construct a further \emph{second-order} identity type of form  $$x_{1}=_{x=_{X}y}y_{1}$$ that says that paths $x_{1},y_{1}$ are, in fact, the same.  The corresponding proof terms are $$x_{2},y_{2} : x_{1}=_{x=_{X}y}y_{1}$$. Under the homotopical interpretation the second-order identity type is interpreted as the space of paths between points $x,y$, while terms $x_{2},y_{2}$ are interpreted as \emph{homotopies} between paths $x_{1},y_{1}$, i.e. as ``paths between the paths'', which intuitively can be thought of as \emph{surfaces} subtended on paths $x_{1},y_{1}$. The process of syntactic building higher identity types can be continued indefinitely. While the intended interpretation of MLTT leaves the semantics of this structure unclear, the homotopical interpretation interprets it as the structure of higher-order fundamental groupoid (which in the general case is an $\infty$-groupoid) of the underlying space $X$. 

The homotopical interpretation of higher identity types allowed Voevodsky to classify types in MLTT using the following inductive definition:

\underline{Definition:}  $S$ is a space (homotopy type) of $h$-level (for ``homotopy level'') $n+1$ if for all its points (terms) $x, y$ path spaces (identity types) $x \ =_{S}\ y$ are of $h$-level $n$

Now we set the $h$-level of \emph{point} (= contractible space) equal to (-2) (to accord with the usual ) and obtain the following stratification of spaces (homotopy types):

\begin{itemize}
\item $h$-level (-2): single point  $pt$;
\item  $h$-level (-1): the empty space $\emptyset$ and the point $pt$: truth-values aka (mere) propositions 
\item   $h$-level 0: sets (discrete point spaces)
\item  $h$-level 1:  flat path groupoids :  no non-contractibe surfaces
\item   $h$-level 2: 2-groupoids : paths and surfaces but no non-contractible volumes
\item
\item
\item  $h$-level $n$: $n$-groupoids
\item $\dots$
\item  $h$-level $\omega$: $\omega$-groupoids
\end{itemize}

Remark that the obtained hierarchy of types is cumulative in the sense that all types of $h$-level $n$ also qualify as types of level $m$ for all $m>n$ (for example, the empty set qualifies as a set, as empty groupoid, etc). So a type of  $h$-level $n$ (aka a $n$-type) can be described as one where the structure of identity types of its point is non-trivial up to level $n$ and is trivial at all higher dimensions.  Given $n$-type $X$ and $k<m$ one can consider $k$-type $\Arrowvert X \Arrowvert_{k}$ obtained via the operation of \emph{truncation}, which amounts to colliding all higher-order terms up to the fixed lower level $k$. In particular, a \emph{propositional truncation} of $n$-type $X$ (where $n> -1$) amounts to determining whether this type is empty or non-empty.    

This hierarchy suggests a modification of the intended interpretation of MLTT along the following lines: instead of interpreting \emph{every} given type $X$ alternatively (or simultaneously as does Martin-L\"of in \cite[p. 13]{Martin-Lof:1984}) as a proposition or as a set, we now 

\begin{itemize}
\item
interpret as \emph{propositions} only types of $h$-level (-1) that have at most one term;
\item
interpret as \emph{sets} only types of $h$-level 0 that have no non-trivial path spaces (i.e. up to homotopy have at most one path between their points)\footnote{The relevant notion of set  is suggested by the identification of sets with discrete topological spaces, and so it differs from the concept of set in ZF and akin set theories. };

\item for $n>0$, interpret $n$-types as fundamental $n$-groupoids.  
\end{itemize}

This new interpretation agrees with Martin-L\"of's intended interpretation of MLTT at the propositional level: a term of proposition (i.e., of (-1)-type) is interpreted as its proof (which, if exists, is unique) as before.  In this case one is, once again, in a position to explain a given proposition $P$ as an intention to give its proof or, equivalently, as the problem ``to prove $P$''. But at the higher $h$-levels the two interpretations diverge:  a set (i.e., a 0-type) cannot be homotopically interpreted as a proposition unless it is a singleton or the empty set, a 1-type (a flat groupoid) cannot be interpreted as a set unless it is a trivial groupoid (with no non-trivial paths between its points), etc. \textemdash\ albeit such reductions to lower-level types are available via an appropriate truncation. Now remark that being a constructive theory, HoTT gives good reasons to think of \emph{every} type (and not only of the propositional types) as a \emph{problem}, and think of terms of this type as  \emph{constructions} that solve it.  As I have already informally explained it above, a construction of some term of 0-type (i.e. of an element of a \emph{set}) under the homotopical interpretation (which, recall, preserves  Martin-L\"of's 1980 interpretation for propositions) requires exhibiting for all terms of this type the relevant identity conditions (in terms of existence of paths). In case of 1-types (flat groupoids) one has a further task to show how different identity proofs (i.e., different paths between its points) are identified and distinguished. This is done in terms of existence of homotopies between these paths. One proceeds similarly in case of higher types using higher homotopies. For a more detailed account of higher-order constructions in HoTT see \cite{Buchholz:2019}. 

Now we get the following picture. A proposition $P$ is identified with the problem ``to prove $P$'' in full accordance with Heyting and Martin-L\"of's \cite{Martin-Lof:1984}. But it is not the case in the given setting that every problem is reducible to this propositional form; higher-order problems that require to build higher-order constructions (starting with set-level constructions) do not reduce in this way (unless one applies the propositional truncation that simplifies the problem). Thus HoTT/UF setting provides a new and independent justification for Kolmogorov's view according to which problems and propositions should not be conflated. The homotopical interpretation of MLTT shows that Kolmogorov is right, and that the distinction between problems and more general propositions is logical rather than linguistic (as far as we qualify HoTT as a logical machinery). The idea that only some but not all problems have form ``to prove proposition $P$'' is compatible with the content of Kolmogorov's 1932 paper \cite{Kolmogoroff:1932} and with the related sources analysed in the present paper. Let me stress once again, that the HoTT-based interpretation of Kolmogorov's distinction between (general) problems and proposition remains constructive and does not use the idea that the logic of propositions unlike the logic of problems is classical.

\section{Conclusion}  
I hope to have shown in this paper that Kolmogorov's 1932 paper \cite{Kolmogoroff:1932}  is an original and seminal work, which is not fully taken into account via the popular BHK-interpretation of intuitionistic logic. The notion of BHK-interpretation is fully justified theoretically as a useful combination of ideas of different thinkers but this combination should not be seen as a final synthesis, which comprises all valuable contents of Kolmogorov's Calculus of Problems. Neither the BHK-interpretation is a reliable historical guide for tracing the historical development of Kolmogorov's logical ideas. In the real history these ideas have an independent life and independent legacy and significance (see  \textbf{4.1.}.   

There are important differences between Kolmogorov's and Heyting's takes on the intuitionistic logic, which I tried to emphasise by talking about the \emph{Kolmogorov-Heyting Controversy} (see  \textbf{3.2.}). This controversy reflects the fruitful discussion between the two thinkers in early 1930s. The controversy is both mathematical and philosophical. Philosophically, Kolmogorov is not sympathetic to Brouwer's and Heyting's idea to rebuild all of mathematics on new intuitionistic principles \textemdash\ whether one thinks about the resulting ``intuitionistic mathematics'' as the only sort of good mathematics, or as a special sort of mathematics co-existing with mathematics of different sorts. Instead, Kolmogorov is looking for applications of Brouwer's and Heyting's logical ideas in his contemporary mathematics as it is practiced without trying to identify its limited fragment that could be called ``intuitionistic''. Mathematically and logically, the divergence between Heyting and Kolmogorov concerns the question of distinguishing between problems and propositions (theorems): while in Heyting's mature view the distinction is not logical but rather merely linguistic, Kolmogorov insists that it is logically significant. 

In the second part of the paper I reviewed a recent work by Sergei Melikhov who justifies Kolmogorov's distinction by developing an original combined logical calculus of problems and propositions (see \textbf{4.1.}); then I proposed another anachronistic justification of the same distinction in terms of Homotopy Type theory (see \textbf{4.2.}). I didn't try to read into Kolmogorov's 1932 paper a pointer to the Homotopy theory; there is no mentioning of homotopy either in this paper or in related Kolmogorov's texts. But the fact that in his philosophical writings Kolmogorov was thinking hard about the concept of mathematical \emph{object} (see \textbf{2.5.}) makes me to believe that he would be deeply interested in HoTT and Univalent Foundations because these recent theories shed a new light on this traditional concept.  It remains an interesting open question whether Melikhov's  combined logic of problems of propositions has some theoretical connections with the proposed HoTT-based justification of Kolmogorov's view on problems and propositions, or these modern developments of Kolmogorov's logical ideas are rather somewhat orthogonal.

\bigskip
\subsection*{Acknowledgements}. I thank Sergei Artemov, Mark van Atten, Lev Beklemishev, Walter Dean, Slava Gerovitch, Lev Lamberov, Sergei Melikhov, Philippe Nabonnand, George Shabat and Noson Yanofsky for their very valuable criticisms and suggestions.



\bibliographystyle{plain} 
\bibliography{kolmoeng} 

\begin{thebibliography}{10}

\bibitem{Aczel:2001}
P.~Aczel.
\newblock The {R}ussell-{P}rawitz modality,.
\newblock {\em Math. Structures Comput. Sci.}, 11:541--554, 2001.

\bibitem{Ahrens&North:2019}
B.~Ahrens and P.~North.
\newblock Univalent {F}oundations and {E}quivalence {P}rinciple.
\newblock {\em S. Centrone, D. Kant and D. Sarikaya (eds.) Reflections on the
  Foundations of Mathematics: Univalent Foundations, Set Theory and General
  Thoughts, Springer, Synthese Library vol.407}, pages 137--150, 2019.

\bibitem{Artemov&Fitting:2019}
S.~Artemov and M.~Fitting.
\newblock {\em Justification {L}ogic: {R}easoning with {R}easons}.
\newblock Cambridge University Press, 2019.

\bibitem{Artemov&Protopopescu:2014}
S.~Artemov and T.~Protopopescu.
\newblock {\em Intuitionistic {E}pistemic {L}ogic}.
\newblock https://arxiv.org/abs/1406.1582v2, 2014.

\bibitem{Artemov:1994a}
S.N. Artemov.
\newblock Logic of {P}roofs.
\newblock {\em Annals of Pure and Applied Logic}, 67:25--59, 1994.

\bibitem{Artemov:2004}
S.N. Artemov.
\newblock Kolmogorov and {G}\''odel’s approach to intuitionistic logic:
  current developments.
\newblock {\em Russian Mathematical Surveys}, 59(2):203--229, 2004.

\bibitem{Artemov:2008}
S.N. Artemov.
\newblock The {L}ogic of {J}ustification.
\newblock {\em The Review of Symbolic Logic}, 4(4):477--513, 2008.

\bibitem{Artemov&Beklemeshev:2004}
S.N. Artemov and L.D. Beklemeshev.
\newblock Provability {L}ogic.
\newblock {\em in: D. Gabbay and F. Guenthner, (eds.) Handbook of Philosophical
  Logic, Second Edition, Dordrecht: Kluwer}, 13:229--403, 2004.

\bibitem{Balaguer:2018}
M.~Balaguer.
\newblock Fictionalism in the {P}hilosophy of {M}athematics.
\newblock {\em Stanford Encyclopedia of Philosophy
  (https://plato.stanford.edu/entries/fictionalism-mathematics/)}, 2018.

\bibitem{Basu&Simpson:2016}
S.S. Basu and S.G. Simpson.
\newblock Mass {P}roblems and {I}ntuitionistic {H}igher-{O}rder {L}ogic.
\newblock {\em Computability}, 5(1):29--47, 2016.

\bibitem{Benacerraf:1965}
P.~Benacerraf.
\newblock What {N}umbers {C}ould {N}ot {B}e.
\newblock {\em Philosophical Review}, 74:47--73, 1965.

\bibitem{Brouwer:1975}
L.E.J. Brouwer.
\newblock {\em Collected {W}orks (ed. by {A.} {H}eyting), vol. 1}.
\newblock North Holland, 1975.

\bibitem{Brouwer:1998a}
L.E.J. Brouwer.
\newblock Does {E}very {R}eal {N}umber {H}ave a {D}ecimal {E}xpansion?
\newblock {\em P. Mancosu (ed.) From Brouwer to Hilbert. The Debate of the
  Foundations of Mathematics in the 1920s. Oxford University Press}, pages
  28--35, 1998.

\bibitem{Brouwer:1998}
L.E.J. Brouwer.
\newblock Mathematica, {S}cience, and {L}anguage.
\newblock {\em P. Mancosu (ed.) From Brouwer to Hilbert. The Debate of the
  Foundations of Mathematics in the 1920s. Oxford University Press}, pages
  45--53, 1998.

\bibitem{Buchholz:2019}
U.~Buchholz.
\newblock Higher-{O}rder {S}tructures in {H}omotopy {T}ype {T}heory.
\newblock {\em S. Centrone, D. Kant and D. Sarikaya (eds.) Reflections on the
  Foundations of Mathematics: Univalent Foundations, Set Theory and General
  Thoughts, Springer, Synthese Library vol.407}, pages 151--172, 2019.

\bibitem{Coquand:2007}
Th. Coquand.
\newblock Kolmogorov's {C}ontribution to {I}ntuitionistic {L}ogic.
\newblock {\em in: E. Charpentier, A. Lesne and N.K. Nikolski (eds).
  Kolmogorov's Heritage in Mathematics}, pages 19--40, 2007.

\bibitem{Fairtlough&Walton:1997}
M.~Fairtlough and M.~Walton.
\newblock Quantified lax logic.
\newblock {\em Tech. report CS-97-11, Univ. of Sheffield, Dept. of Computer
  Science}, 1997.

\bibitem{Frege:1903}
G.~Frege.
\newblock {\em Grundgesetze der {A}rithmetik, {B}and 2}.
\newblock Olms, 1962.

\bibitem{Friedberg:1956}
R.M. Friedberg.
\newblock The {S}olution of {P}ost's {P}roblem.
\newblock {\em Bulletin of American Mathematical Society}, 62(3):49--59, 1956.

\bibitem{Godel:1933}
K.~G\"odel.
\newblock Zur intuitionistischen {A}rithmetik und {Z}ahlentheorie.
\newblock {\em Ergebnisse eines Math. Kolloquimus}, 4:34--38, 1933.

\bibitem{Godel:1933a}
K.~G\"odel.
\newblock Eine {I}nterpretation des intuitionistischen {A}ussagenkalkuls.
\newblock {\em Ergebnisse eines Math. Kolloquimus}, 14:39--40, 1933a.

\bibitem{Godel:1986}
K.~Godel.
\newblock {\em Collected {W}orks, ed. {F}eferman et al., Vol. 1, {P}ublications
  1929-1936}.
\newblock Oxford University Press, 1986.

\bibitem{Grayson:2018}
D.~Grayson.
\newblock An {I}ntroduction of {U}nivalent {F}oundations for {M}athematicians.
\newblock {\em Bulletin of American Mathematical Society (New Series)},
  55(4):427--450, 2018.

\bibitem{UF:2013}
Univalent~Foundations Group.
\newblock {\em Homotopy {T}ype {T}heory: {U}nivalent {F}oundations of
  {M}athematics}.
\newblock Institute for Advanced Study (Princeton); available at
  http://homotopytypetheory.org/book/, 2013.

\bibitem{Heyting:1930}
A.~Heyting.
\newblock Die formalen {R}egeln der intuitionistischen {L}ogik {I}.
\newblock {\em Sitzungsberichte der Preussischen Akademie der Wissenschaften},
  pages 42--56, 1930.

\bibitem{Heyting:1930a}
A.~Heyting.
\newblock Sur la logique intuitionniste.
\newblock {\em Acad\'emie Royale de Belgique, Bulletin de la Classe des
  Sciences}, 16:957--963, 1930.

\bibitem{Heyting:1931}
A.~Heyting.
\newblock Die intuitionistische {G}rundlegung der {M}athematik.
\newblock {\em Erkenntnis}, 2:106--115, 1931.

\bibitem{Heyting:1934}
A.~Heyting.
\newblock {\em Mathematische {G}rundlagenforschung, {I}ntuitionismus,
  {B}eweistheorie}.
\newblock Springer, 1934.

\bibitem{Heyting:1956}
A.~Heyting.
\newblock {\em Intuitionism: and introduction}.
\newblock North Holland, 1956.

\bibitem{Heyting:1958}
A.~Heyting.
\newblock Intuitionism in {M}athematics.
\newblock {\em R. Klibansky (ed.), La philosophie au milieu du vingti\`eme
  si\`ecle, Firenze: La nuova Italia}, pages 101--115, 1958.

\bibitem{Heyting:1983}
A.~Heyting.
\newblock The intuitionist foundations of mathematics.
\newblock {\em P. Benaceraff and H. Putnam (eds.), Philosophy of Mathematics.
  Selected Readings. Cambridge University Press}, pages 52--60, 1983.

\bibitem{Heyting:1998b}
A.~Heyting.
\newblock On {F}ormal {R}ules of {I}ntuitionistic {L}ogic.
\newblock {\em P. Mancosu (ed.) From Brouwer to Hilbert. The Debate of the
  Foundations of Mathematics in the 1920s. Oxford University Press}, pages
  311--327, 1998.

\bibitem{Heyting:1998a}
A.~Heyting.
\newblock On {I}ntuitionistic {L}ogic.
\newblock {\em P. Mancosu (ed.) From Brouwer to Hilbert. The Debate of the
  Foundations of Mathematics in the 1920s. Oxford University Press}, pages
  306--310, 1998.

\bibitem{Heyting:1998}
A.~Heyting.
\newblock On the {F}ormal {R}ules of {I}ntuitionistic {L}ogic.
\newblock {\em P. Mancosu (ed.) From Brouwer to Hilbert. The Debate of the
  Foundations of Mathematics in the 1920s. Oxford University Press}, pages
  311--327, 1998.

\bibitem{Hilbert:1899}
D.~Hilbert.
\newblock {\em Grundlagen der {G}eometrie}.
\newblock Leipzig, 1899.

\bibitem{Hilbert:1928}
D.~Hilbert.
\newblock Die {G}rundlagen der {M}athematik.
\newblock {\em Abhandlungen aus dem Seminar der Hamburgischen Universit\"at},
  6:65--85, 1928.

\bibitem{Hilbert:1967}
D.~Hilbert.
\newblock Foundations of {M}athematics.
\newblock {\em J. van Heijenoort (ed.), From Frege to G\"odel: A Source Book in
  the Mathematical Logic}, 2:464--480, 1967.

\bibitem{Hilbert&Bernays:1934-1939}
D.~Hilbert and P.~Bernays.
\newblock {\em Grundlagen der {M}athematik}.
\newblock Springer, 1934-1939.

\bibitem{Jacquette:1996}
D.~Jacquette.
\newblock On defoliating meinong's jungle.
\newblock {\em Axiomathes}, 7:17--42, 1996.

\bibitem{Japaridze:2003}
G.~Japaridze.
\newblock Introduction to computability logic.
\newblock {\em Annals of Pure and Applied Logic}, 123:1--99, 2003.

\bibitem{Japaridze:2007}
G.~Japaridze.
\newblock Intuitionistic computability logic.
\newblock {\em Acta Cybernetica}, 18(1):77--113, 2007.

\bibitem{Khinchin:1926}
A.Ya. Khinchin.
\newblock Idei intuitsionizma i bor'ba za predmet v sovremennoi matematike
  [ideas of intuitionism and the struggle over subject matter in contemporary
  mathematics]({R}ussian).
\newblock {\em Vestnik Kommunisticheskoi akademii [Messenger of Communist
  Academy]}, pages 184--192, 1926.

\bibitem{Kleene:1945}
S.C. Kleene.
\newblock On the {I}nterpretation of {I}ntuitionistic {N}umber {T}heory.
\newblock {\em The Journal of Symbolic Logic}, 10(4):109--124, 1945.

\bibitem{Kleene:1973}
S.C. Kleene.
\newblock Realizability: a {R}etrospective {S}urvey.
\newblock {\em in: A.R.D. Mathias and H. Rogers (eds.) Cambridge Summer School
  in Mathemamtical Logic 1971 (Lecture Notes in Mathematics)}, pages 95--112,
  1973.

\bibitem{Kolmogoroff:1932}
A.~Kolmogoroff.
\newblock Zur {D}eutung der {I}ntuitionistischen {L}ogik.
\newblock {\em Mathematische Zeitschrift}, 35:58--65, 1932.

\bibitem{Kolmogoroff:1937}
A.~Kolmogoroff.
\newblock Review of ``zur intuitionistische {D}eutung logischer {F}ormeln'' by
  {H.} {F}redeuntal and ``{B}emerkungen zu dem {A}ugsatz von {H}errn
  {F}redeuntal '{Z}ur intuitionistische {D}eutung logischer {F}ormeln' by {A}.
  {H}eyting''.
\newblock {\em Zentralblatt f\"ur Mathematik und ihre Grenzgebiete}, page
  0015.24201, 1937.

\bibitem{Kolmogorov:1925}
A.N. Kolmogorov.
\newblock O principe tertium non datur ({R}ussian).
\newblock {\em Matematiceskij Sbornik}, 32:646--667, 1925.

\bibitem{Kolmogorov:1929}
A.N. Kolmogorov.
\newblock Contemporary {D}ebates on the {N}ature of {M}athematics ({R}ussian).
\newblock {\em Novoe Slovo}, 6:41--54, 1929.

\bibitem{Kolmogorov:1936}
A.N. Kolmogorov.
\newblock Preface ({R}ussian).
\newblock {\em Russian edition of 'Mathematische Grundlagenforschung,
  Intuitionismus, Beweistheorie' by A. Heyting, Springer 1934}, pages 3--4,
  1936.

\bibitem{Kolmogorov:1988}
A.N. Kolmogorov.
\newblock Letters of {A}.{N}. {K}olmogorov to {A}. {H}eyting, transl. and comm.
  {V}.{E}. {P}lisko.
\newblock {\em Russian Math. Surveys}, 43(6):89--93, 1988.

\bibitem{Kolmogorov:1991}
A.N. Kolmogorov.
\newblock On the papers on {I}ntuitionistic {L}ogic.
\newblock {\em V.M. Tikhomirov (ed.) Selected Works of A.N. Kolmogorov,
  Springer}, 1:451--452, 1991.

\bibitem{Kolmogorov:1991a}
A.N. Kolmogorov.
\newblock On the tertium non datur principle.
\newblock {\em V.M. Tikhomirov (ed.) Selected Works of A.N. Kolmogorov,
  Springer}, 1:40--68, 1991.

\bibitem{Kolmogorov:1998}
A.N. Kolmogorov.
\newblock On the {I}nterpretation of {I}ntuitionistic {L}ogic.
\newblock {\em P. Mancosu (ed.) From Brouwer to Hilbert. The Debate of the
  Foundations of Mathematics in the 1920s. Oxford University Press}, pages
  328--334, 1998.

\bibitem{Kolmogorov:2006}
A.N. Kolmogorov.
\newblock Contemporary {D}ebates on the {N}ature of {M}athematics.
\newblock {\em Problems of Information Transmission}, 42(4):379--389, 2006.

\bibitem{Kozen:2006}
D.C. Kozen.
\newblock {\em Theory of {C}omputation}.
\newblock Springer, 2006.

\bibitem{Martin-Lof:1984}
P.~Martin-L\"of.
\newblock {\em Intuitionistic {T}ype {T}heory (Notes by {G}iovanni {S}ambin of
  a series of lectures given in {P}adua, {J}une 1980)}.
\newblock Napoli: BIBLIOPOLIS, 1984.

\bibitem{Medvedev:1955}
Yu.~T. Medvedev.
\newblock Degrees of difficulty of the mass problems.
\newblock {\em Doklady Academii Nauk SSSR (Russian)}, 104(4):501--504, 1955.

\bibitem{Medvedev:1962}
Yu.~T. Medvedev.
\newblock Finitary {P}roblems.
\newblock {\em Doklady Academii Nauk SSSR (Russian)}, 142(5):1015--1018, 1962.

\bibitem{Medvedev:1963}
Yu.~T. Medvedev.
\newblock Interpretation of {L}ogical {F}ormulas with {F}initary {P}roblems and
  {I}ts {R}elation to the {R}ealisability {T}heory.
\newblock {\em Doklady Academii Nauk SSSR (Russian)}, 148(4):771--774, 1963.

\bibitem{Medvedev:1966}
Yu.~T. Medvedev.
\newblock On the {I}nterpretation of {L}ogical {F}ormulas with {F}initary
  {P}roblems.
\newblock {\em Doklady Academii Nauk SSSR (Russian)}, 169(1):20--23, 1966.

\bibitem{Meinong:1904}
A.~Meinong.
\newblock \"uber {G}egenstandstheorie.
\newblock {\em in A. Meinong (Ed.) Untersuchungen zur Gegenstandstheorie und
  Psychologie, Leipzig, Verlag von Johann Ambrosius Barth 1904}, pages 1--50,
  1904.

\bibitem{Meinong:1907}
A.~Meinong.
\newblock {\em \"Uber die {S}tellung der {G}egenstandstheorie im {S}ystem der
  {W}issenschaften}.
\newblock Lepzig: R. Voigtl\"ander, 1907.

\bibitem{Meinong:1960}
A.~Meinong.
\newblock The {T}heory of {O}bjects.
\newblock {\em in R. Chisholm, Realism and the Background of Phenomenology,
  Free Press of Glencoe, Illinois, 1960}, pages 76--117, 1960.

\bibitem{Melikhov:2017}
S.A. Melikhov.
\newblock {\em Mathematical {S}emantics of {I}ntuitionistic {L}ogic}.
\newblock https://arxiv.org/abs/1504.03380, 2017.

\bibitem{Melikhov:2022a}
S.A. Melikhov.
\newblock {\em A {G}alois connection between classical and intuitionistic
  logics. I: {S}yntax}.
\newblock https://arxiv.org/abs/1504.03380, 2022.

\bibitem{Melikhov:2022b}
S.A. Melikhov.
\newblock {\em A {G}alois connection between classical and intuitionistic
  logics. II: {S}emantics}.
\newblock https://arxiv.org/abs/1504.03380, 2022.

\bibitem{Melikhov:2023}
S.A. Melikhov.
\newblock A {J}oint {L}ogic of {P}roblems and {P}ropositions.
\newblock {\em Doklady Mathematics, in press}, 2023.

\bibitem{Mostowski:1956}
A.~Mostowski.
\newblock Review of '{D}egrees of difficulty of the mass problem' by {Y}u.{T}.
  {M}edvedev,.
\newblock {\em Journal of Symbolic Logic}, 21(3):320--321, 1956.

\bibitem{Muchnik:1956}
A.A. Muchnik.
\newblock The insolubility of the reduction problem in the theory of
  algorithms.
\newblock {\em Doklady Academii Nauk SSSR (Russian)}, 108(2):194--197, 1956.

\bibitem{Muchnik:1959}
A.A. Muchnik.
\newblock The {P}roblem of {R}educibility for {R}ecursive sets.
\newblock {\em Mathematical Enlightment (Russian)}, 4:233--236, 1959.

\bibitem{Muchnik:1963}
A.A. Muchnik.
\newblock Strong and weak reducibility of algorithmic problems.
\newblock {\em Siberian Mathematical Journal (Russian)}, 4:1328--1341, 1963.

\bibitem{Muchnik:2016}
A.A. Muchnik.
\newblock Strong and weak reducibility of algorithmic problems.
\newblock {\em Computability}, 5(1):49--59, 2016.

\bibitem{Onoprienko:2020}
A.A. Onoprienko.
\newblock Kripke-type semantics for a logic of problems and propositions.
\newblock {\em Sbornik Mathematics}, 211:709--732, 2020.

\bibitem{Onoprienko:2022}
A.A. Onoprienko.
\newblock Topological models of propositional logic of problems and
  propositions.
\newblock {\em Moscow University Mathematics Bulletin}, 77(5):236--241, 2022.

\bibitem{Rodin:2014}
A.~Rodin.
\newblock {\em Axiomatic {M}ethod and {C}ategory {T}heory ({S}ynthese {L}ibrary
  vol. 364)}.
\newblock Springer, 2014.

\bibitem{Routley&Routley:1973}
R.~Routley and V.~Routley.
\newblock Rehabilitating {M}einong's {T}heory of {O}bjects.
\newblock {\em Revue Internationale de Philosophie}, 27(104/105):224--254,
  1973.

\bibitem{Simpson:2008}
S.G. Simpson.
\newblock Weak {C}ounter-{E}xamples.
\newblock {\em Notre Dame Journal of of Formal Logic}, 49:127--136, 2008.

\bibitem{Skvortsov:1979a}
D.P. Skvortsov.
\newblock Logic of {I}nfinitary {P}roblems and {K}ripke {M}odels on {A}tomic
  {S}emi-{L}attices of {S}ets.
\newblock {\em Doklady Academii Nauk SSSR (Russian)}, 245(4):798--801, 1979.

\bibitem{Skvortsov:1979}
D.P. Skvortsov.
\newblock Two {G}eneralisations of the {C}oncept of {F}initary {P}roblem.
\newblock {\em in: A.I. Mikhailov (ed.) Studies in Non-Classical Logics and Set
  Theory, Nauka Publishing, 1979 (Russian)}, pages 201--240, 1979.

\bibitem{Stone:1961}
M.~Stone.
\newblock Revolution in {M}athematics.
\newblock {\em American Mathematical Monthly}, 68(8):715--734, 1961.

\bibitem{Sundholm:1983}
G.~Sundholm.
\newblock Constructions, {P}roofs, and the {M}eaning of {L}ogical {C}onstants.
\newblock {\em Journal of Philosophical Logic}, 12(2):151--172, 1983.

\bibitem{Sundholm:2012}
G.~Sundholm.
\newblock On the {P}hilosophical {W}ork of {P}er {M}artin-{L}\"of.
\newblock {\em Preface to: Dybjer P. et al. (Eds.), Epistemology versus
  Ontology Essays on the Philosophy and Foundations of Mathematics in Honour of
  Per Martin-L\"of. Logic, Epistemology, and the Unity of Science no. 27.
  Dordrecht: Springer Netherlands}, pages xvii--xxiv, 2012.

\bibitem{Troelstra:1990}
A.S. Troelstra.
\newblock On the {E}arly {H}istory of {I}ntuitionistic {L}ogic.
\newblock {\em P.P. Petkov (ed.) Mathematical Logic, New York: Plenum Press},
  pages 3--17, 1990.

\bibitem{Uspenskii:2006}
V.A. Uspenskii.
\newblock Kolmogorov as {I} {R}emember {H}im (russian).
\newblock {\em A.N. Shiryaev (ed.), Kolmogorov in the Memories of his Pupils
  (Russian), Moscow Centre of Continuous Mathematical Education}, pages
  272--371, 2006.

\bibitem{Uspenskii&Plisko:1991}
V.A. Uspenskii and V.E. Plisko.
\newblock Intuitionistic {L}ogic.
\newblock {\em V.M. Tikhomirov (ed.) Selected Works of A.N. Kolmogorov,
  Springer}, 1:452--465, 1991.

\bibitem{vanAtten:2020}
M.~van Atten.
\newblock Weak {C}ounter-{E}xamples.
\newblock {\em Stanford Encyclopedia of Philosophy
  (https://plato.stanford.edu/entries/brouwer/weakcounterex.html), Supplement
  to the main entry 'Luitzen Egbertus Jan Brouwer',
  https://plato.stanford.edu/entries/brouwer/index.html}, 2020.

\bibitem{vanAtten:2022}
M.~van Atten.
\newblock The {D}evelopment of {I}ntuitionistic {L}ogic.
\newblock {\em Stanford Encyclopedia of Philosophy
  (https://plato.stanford.edu/entries/intuitionistic-logic-development/)},
  2022.

\bibitem{vanStigt:1990}
W.P. van Stigt.
\newblock {\em Brouwer's {I}ntuitionism}.
\newblock North Holland, 1990.

\bibitem{Verburgt&Hoppe-Kondrikova:2006}
L.M. Verburgt and O.~Hoppe-Kondrikova.
\newblock On {A}.{Y}a. {K}hinchin’s paper ‘ideas of intuitionism and the
  struggle for a subject matter in contemporary mathematics’ (1926): A
  translation with introduction and commentary.
\newblock {\em Historia Mathematica}, 43(4):369--398, 2006.

\bibitem{Voevodsky:2006}
V.A. Voevodsky.
\newblock Foundations of {M}athematics and {H}omotopy {T}heory (lecture at the
  {P}rinceton institute of {A}dvanced {S}tudies on march 22, 2006), 2006.
\newblock slides available online via
  https://www.math.ias.edu/vladimir/Lectures.

\bibitem{Voevodsky:2010a}
V.A. Voevodsky.
\newblock Univalent {F}oundations (lecture at the {P}rinceton institute of
  {A}dvanced {S}tudies on {D}ecember 10, 2010), 2010.
\newblock video available online at
  https://www.ias.edu/video/univalent/voevodsky.

\bibitem{Voevodsky:2010}
V.A. Voevodsky.
\newblock Univalent {F}oundations (lecture in {B}onn on {S}eptember 8, 2010),
  2010.
\newblock slides available online via
  https://www.math.ias.edu/vladimir/Lectures.

\end{thebibliography}

\clearpage

\section*{Appendix:  A.N. Kolmogorov, \emph{Preface} to  Russian edition  \cite{Kolmogorov:1936} of \emph{Mathematische Grundlagenforschung, Intuitionismus, Beweistheorie} by Arend Heyting  \cite{Heyting:1934}. Translation from Russian by Andrei Rodin.}

\bigskip
 In the Introduction to his book Heyting quite rightly states that the purpose of research in the foundations of mathematics is not limited to the verification of their firmness, the abolishing  of ill-founded parts of mathematics and the elimination of emerging contradictions but also includes a positive analysis of the subject-matter of mathematics, of its methods and ways of its development. Nevertheless, since the efforts of mathematicians researching the foundations of their science during the last years were focused mostly on the former goal, Heyting in the present review pays little attention to the second goal. Thus one should keep firmly in mind
that Heyting’s book is not an introduction into a positive philosophy of mathematics.  It is, by and large, an exposition of various directions and ways of critique concerning some [theoretical] constructions used by working mathematicians, and of attempts to mitigate the destructive results of this critique. Heyting does not limit his exposition to an analysis of the two major trends in this area of research (Intuitionism and Formalism) but also makes interesting remarks pointing to their possible synthesis.

The destructive, critical parts of both these trends are essentially related to the principles of Subjective Idealism. This is why Heyting’s review cannot give any indication to how the result of critique of basic mathematical concepts would look like if the critique proceeds on the materialistic grounds.  But the constructive side of the two trends (development of mathematics without using the Law of Excluded Middle in the case of Intuitionism, and methods of proving consistency in the case of Formalism) comes down to a concrete mathematical work, which is often admirably sharp and smart. In spite of their [erroneous] philosophical subjective idealistic assumptions,  the two schools [of mathematical thought] discovered in the course of this work  a number of extraordinarily deep and interesting facts.

We cannot agree with the intuitionists when they claim that mathematical objects are products of the constructive activity of our spirit. For us, mathematical objects are abstractions from existing forms of reality, which is independent from our spirit. We know that the constructive solutions of problems are as much important in mathematics as the pure proofs of theoretical sentences. This constructive aspect of mathematics does not conceal for us its other and more fundamental aspect, namely, its epistemic aspect.  But the laws of mathematical construction discovered by Brouwer and systematized by Heyting under the appearance of new intuitionistic logic, so understood, preserve for us their fundamental significance (see my article \cite{Kolmogoroff:1932}).  

We equally reject the tenets of the formalistic school. We believe that mathematics does not reduce to its formalised part, and that [mathematical] problems that cannot be solved on the basis of axioms formulated so far, [may] nevertheless have a well-determined real sense and admit unambiguous solutions. This view, however, does not diminish our interest to studying the structure of formalised mathematical theories. 

Every \emph{accomplished} mathematical theory should be formalised because the final goal of every such theory is building an algorithm for an automatic solution of the related problems. It is moreover interesting that the studies made by the formalists themselves  lead one to the conclusion that, first, the formalisation [of mathematics] can never be completed and, second, that the formalisation needs some contentual (not formal) mathematics (under the name of “metamathematics”, if one wishes) as its prerequisite. 

A danger of using Heyting’s review lies in its extreme density. This density makes it neccessary for the reader to consult the primary sources listed in the Bibliography, which is carefully composed and complete  ([for publications made] before 1933). This density also makes the task of translator [of Heyting’s book into Russian] very difficult. This difficult task was accomplished by the translator (A.P. Yushkevich) perfectly. He always aimed at the precise translation of the author’s thoughts, even at the price of some additional length. When this is needed, the newly introduced Russian equivalents of some terms are accompanied by the original [German] words in parentheses.

\end{document}